\documentclass{amsart}

\usepackage{latexsym}
\usepackage{palatino}

\newcommand{\SP}{\operatorname{Sp}}
\newcommand{\GL}{\operatorname{GL}}
\newcommand{\Aut}{\operatorname{Aut}}

\newcommand{\SL}{\operatorname{SL}}

\newcommand{\Hom}{\operatorname{Hom}}

\newcommand{\End}{\operatorname{End}}
\newcommand{\tensor}{\otimes}

\newcommand{\coker}{\operatorname{coker}}
\newcommand{\imag}{\operatorname{Im}}
\newcommand{\Dist}{\operatorname{Dist}}

\DeclareSymbolFont{curgreek}{U}{eur}{m}{n}
\DeclareMathSymbol{\varpi}{\mathord}{curgreek}{"24}

\DeclareMathSymbol{\cursivebeta}{\mathord}{curgreek}{"0C}
\DeclareMathSymbol{\cursiveepsilon}{\mathord}{curgreek}{"22}
\DeclareMathSymbol{\cursiveeta}{\mathord}{curgreek}{"11}
\DeclareMathSymbol{\cursivepsi}{\mathord}{curgreek}{"20}

\DeclareMathAlphabet{\mathrc}{U}{eur}{m}{n}
\DeclareMathAlphabet{\matheucal}{U}{eus}{m}{n}

\newcommand\N{\mathbb{N}}    
\newcommand\Z{\mathbb{Z}}    
\newcommand\F{\mathbb{F}}    
\newcommand\Q{\mathbb{Q}}    
\newcommand\G{\mathbb{G}}    

\newcommand\Aff{\mathbb{A}}      

\newcommand\Witt{\mathrc{W}}   

\newcommand{\mm}{\mathfrak{m}}

\newcommand\Lie{\operatorname{Lie}}
\newcommand\Ad{\operatorname{Ad}}
\newcommand\ad{\operatorname{ad}}

\newcommand\lie[1]{\mathfrak{#1}}
\newcommand\glie{\lie{g}}

\newcommand\plie{\lie{p}}
\newcommand\gl{\lie{gl}}

\newcommand{\hlie}{\lie{h}}
\newcommand{\ulie}{\lie{u}}
\newcommand{\vlie}{\lie{v}}
\newcommand{\alie}{\lie{a}}
\newcommand{\wlie}{\lie{w}}
\newcommand{\sllie}{\lie{sl}}
\newcommand{\clie}{\lie{c}}
\newcommand{\llie}{\lie{l}}
\newcommand{\nlie}{\lie{n}}

\newcommand{\zlie}{\lie{z}}

\newcommand{\normal}{\lhd}

\newcommand{\congruent}{\equiv} 
\newcommand{\iso}{\simeq}

\swapnumbers

\theoremstyle{plain}

\newtheorem*{theorem}{Theorem}

\newtheorem*{prop}{Proposition}

\newtheorem*{cor}{Corollary}

\newtheorem*{lem}{Lemma}

\swapnumbers
\newtheorem{stmt}{\quad}[subsection]

\swapnumbers
\theoremstyle{remark}
\newtheorem*{rem}{Remark}

\newtheorem*{rems}{Remarks}

\newtheorem*{example}{Example}

\numberwithin{equation}{subsection}

\newcommand{\NN}{\mathcal{N}}

\usepackage{amssymb}

\usepackage{xypic}

\begin{document}
\bibliographystyle{amsalpha} 
\author{George J. McNinch}
\thanks{This work was supported by a grant from the 
        National Science Foundation.}
\email{McNinch.1@nd.edu} 
\date{January 25, 2001}

\address{Department of Mathematics \\ 
  University of Notre Dame \\ Notre Dame, IN~ 46556~ USA}

\title{Abelian Unipotent Subgroups of Reductive Groups}

\maketitle

\newcommand{\Levi}{\mathrc{L}}
\newcommand{\V}{\mathrc{V}}
\newcommand{\da}{\bullet}
\newcommand{\A}{\mathcal{A}}
\newcommand{\e}{\varepsilon}
\newcommand{\UU}{\mathcal{U}}
\newcommand{\LL}{\mathcal{L}}
\newcommand{\FF}{\mathrc{F}}
\newcommand{\EE}{\mathrc{E}}
\newcommand{\Zp}{\Z_{(p)}}
\newcommand{\VV}{\mathcal{V}}
\newcommand{\pow}[1]{[ \hspace{-1.6pt} [ {#1} ] \hspace{-1.6pt} ]}
\newcommand{\powfield}[1]{( \hspace{-1.6pt} ( {#1} ) \hspace{-1.6pt} )}

\begin{abstract}
  Let $G$ be a connected reductive group defined over an algebraically
  closed field $k$ of characteristic $p > 0$. The purpose of this
  paper is two-fold. First, when $p$ is a good prime, we give a new
  proof of the ``order formula'' of D. Testerman for unipotent
  elements in $G$; moreover, we show that the same formula determines
  the $p$-nilpotence degree of the corresponding nilpotent elements in the
  Lie algebra $\glie$ of $G$.
  
  Second, if $G$ is semisimple and $p$ is sufficiently large, we show
  that $G$ always has a faithful representation $(\rho,V)$ with the
  property that the exponential of $d\rho(X)$ lies in $\rho(G)$ for
  each $p$-nilpotent $X \in \glie$.  This property permits a
  simplification of the description given by Suslin, Friedlander, and
  Bendel of the (even) cohomology ring for the Frobenius kernels
  $G_d$, $d \ge 2$. The previous authors already observed that the
  natural representation of a classical group has the above property
  (with no restriction on $p$).  Our methods apply to any Chevalley
  group and hence give the result also for quasisimple groups with
  ``exceptional type'' root systems.  The methods give explicit
  sufficient conditions on $p$; for an adjoint semisimple $G$ with
  Coxeter number $h$, the condition $p > 2h -2$ is always good enough.
\end{abstract}

\section{Introduction}

Let $k$ be an algebraically closed field of characteristic $p>0$, and
let $G$ be a connected, reductive group over $k$. We consider in this
paper two questions which involve the relationship between nilpotent
elements in the Lie algebra $\glie$ of $G$ and certain unipotent
subgroups of $G$.

\subsection{}
\label{sub:first_theorem}

There are finitely many (adjoint) orbits of $G$ on the nilpotent
elements of its Lie algebra $\glie$; since $\glie$ is a $p$-Lie
algebra it is reasonable to ask for each nilpotent class $\Ad(G)X
\subset \glie$ what is the minimal integer $m \ge 1$ for which 
$X^{[p^m]} = 0$.

The analogous question for unipotent elements in $G$ was answered in
\cite{testerman}; D. Testerman gave there a formula for the orders of
the unipotent elements in $G$.  We show here that the answer in both
cases is ``the same'' and that moreover by first proving the Lie
algebra result, one obtains a proof of Testerman's formula which avoids
the calculations with explicit representatives for the unipotent
classes that were carried out in \emph{loc. cit.}

More precisely, we prove the following:
\begin{theorem}
  Assume that $p$ is a good prime for the connected reductive group
  $G$, and that $P$ is a distinguished parabolic subgroup of $G$ with
  unipotent radical $\V$. Write $n(P)$ for the nilpotence class of
  $\V$ (which is the same as the nilpotence class of $\vlie$), and let
  the integer $m > 0$ be minimal with the property that $p^m \ge
  n(P)$. 
  \begin{enumerate}
  \item The $p$-nilpotence degree of a Richardson element of $\vlie =
    \Lie(\V)$ is $m$; equivalently, the $p$-exponent of the Lie
    algebra $\vlie$ is $m$;
  \item The order of a Richardson element of $\V$ is $p^m$;
    equivalently, the exponent of $\V$ is $p^m$.
  \end{enumerate}
\end{theorem}

In section \ref{sec:nilpotent} we recall general notions and
definitions concerning nilpotent group and nilpotent Lie algebras.
There is a simple formula for the nilpotence class $n(P)$ given in
section \ref{sub:nilp-class}.  For $m$ as in the theorem, it follows
from generalities (see Lemma \ref{sub:general_nilp_lemma} and
Proposition \ref{sub:nilp-class}) that the $p$-nilpotence degree of
$X$ is $\le m$ and that the order of $u$ is $\le p^m$.  Thus, the
theorem amounts to the following assertions: the exponent of $\V$ and
the $p$-exponent of $\vlie$ are as large as permitted by their
respective nilpotence class.

In section \ref{sec:Abelian-unipotent}, we discuss connected, Abelian,
unipotent algebraic groups. In characteristic 0, any such group is a
vector group, but that is not true in positive characteristic. On the
other hand, in the positive characteristic case, any such group $U$ is
\emph{isogenous} to a product of ``Witt vector groups'' whose
dimensions are uniquely determined by $U$.  Using this we observe that
the $p$-exponent of $\Lie(U)$ is $\le \log_p$ of the
exponent of the group $U$; see Proposition \ref{sub:connect-abelian}.

In section \ref{sec:reductive}, we review relevant facts concerning
the classification of unipotent and nilpotent classes for reductive
groups. The Bala-Carter theorem, as proved for all good primes $p$ by
Pommerening, parameterizes the nilpotent classes in the Lie algebra
$\glie$ of $G$; thanks to a result of Springer, it then also
parameterizes the unipotent classes in $G$.

If $\V$ is as in the theorem, a result of Spaltenstein shows that the
centralizer dimension of a Richardson element $X \in \Lie(\V)$ is
``the same as in characteristic 0''; using this fact, we are able to
use reduction modulo $p$ arguments to obtain a lower bound on the
nilpotence class of $\ad(X)$ which suffices for part (1) of the
theorem; the details are contained in section \ref{sec:Lie-alg}.

Next, we locate a connected, Abelian, unipotent subgroup $Z$ of $G$
which meets the Richardson orbit of $P$ on $\V$, and moreover such
that $\zlie = \Lie(Z)$ meets the Richardson orbit of $P$ on $\vlie$.
The results in section \ref{sec:Abelian-unipotent} show now that
$\log_p$ of the exponent of $Z$ must be $\ge$ the $p$-exponent of
$\zlie$, from which we deduce part (2) of the theorem. This argument
is contained in section \ref{sec:Group}.

\subsection{}
Let $H$ denote a linear algebraic group over $k$ defined over $\F_p$.
In \cite{SFB:infin-1-psg} and \cite{SFB:support-var}, Suslin,
Friedlander, and Bendel relate the cohomology of the Frobenius kernel
$H_d$ to a certain affine scheme $\underline{\A}(d,H)$ whose
$k$-points coincide with the set of all group scheme homomorphisms
$\G_{a,d} \to H$. In fact, they show that the spectrum of the even
cohomology ring of $H_d$ is homeomorphic to $\underline{\A}(d,H)$.
Let $\A(d,H)$ be the variety corresponding to $\underline{\A}(d,H)$
[if $A$ denotes the coordinate ring of the scheme, then the coordinate
ring of the variety is $A'=A/\sqrt{0}$; in this case, the maximal
ideals of $A'$ identify with the above group scheme homomorphisms]. We
observe that $\A(d,H)$ and $\underline{\A}(d,H)$ are homeomorphic, and
in this paper we will only work with the variety.

In case $H$ is the full linear group $\GL(V)$ of a $k$-vectorspace
$V$, $\A(d,H)$ has a simple description as the variety 
of commuting $d$-tuples of $p$-nilpotent elements of $\hlie =
\lie{gl}(V)$.

For general $H$, one may take a faithful representation $(\rho,V)$ of
$H$ and it was observed in \cite{SFB:infin-1-psg} that $\A(d,H)$ is a
somewhat mysterious closed subvariety of $\A(d,\GL(V))$.  If for each
$p$-nilpotent $X \in \glie$ the exponential homomorphism $t \mapsto
\exp(d\rho(tX))$ takes values in $H$, we say that $(\rho,V)$ is an
exponential-type representation. It is shown in \emph{loc.
  cit.}  that if $H$ has an exponential-type representation, then
$\A(d,H)$ may be identified with the variety $\NN_p(d,\hlie)$ of
commuting $d$-tuples of $p$-nilpotent elements in $\hlie$.

We show in this paper that if $G$ is semisimple and $p$ is
sufficiently large, then $G$ has an exponential-type representation.
We consider exponentials in section \ref{sec:exponentials}; the
results on exponentials in Chevalley groups may be found in
\ref{sub:admiss-group}.  As a by-product of some of our constructions,
we obtain also a new proof, for classical groups, of a recent result
of Proud \cite{Proud-Witt} concerning Witt-vector subgroups containing
unipotent elements; see Theorem \ref{sub:classical-witt}.

When $p$ does not divide the order of the ``fundamental group'' of
$G$, we show that $\A(d,G)$ is isomorphic to $\A(d,G_{sc})$ as
$G_{sc}$-varieties, where $G_{sc}$ is the simply connected covering
group; in this sense, $\A(d,G)$ is independent of isogeny. However, it
is not at all clear whether the property of having an exponential-type
representation is independent of isogeny.  

When $G$ is a classical group, it was observed in
\cite{SFB:infin-1-psg} that its ``natural'' module $V$ defines an
exponential-type representation (in any characteristic); so long as $p$
does not divide the order of the fundamental group, this shows that
$\A(d,G') \iso \NN_p(d,\glie)$ for any quasisimple, semisimple group
$G'$ with root system of type $A$, $B$, $C$ or $D$.

For a general semisimple group $G$ we show that if $p > 2h-2$ (where
$h$ is the Coxeter number), the adjoint module is an exponential-type
representation for the corresponding adjoint group (which is isogenous to
$G$); since this inequality also guarantees that $p$ doesn't divide
the order of the fundamental group, we get $\A(d,G) \iso
\NN_p(d,\glie)$ with this condition on $p$.

If $G$ is an exceptional group of type $E_8$, our techniques do no
better than the bound $p> 2h-2 = 58$. For the other exceptional
groups, we improve this bound slightly; see \ref{sub:scheme-iso}.

Suppose that $G$ has an exponential-type representation. As observed
in \cite[Remark 1.9]{SFB:infin-1-psg}, it is not clear whether the
resulting isomorphism $\NN_p(d,\glie) \to \A(d,G)$ is intrinsic, or
depends on the choice of exponential-type representation. In an
attempt to study this question, we consider in section
\ref{sec:exponential-big-p} a related result due to Serre concerning
exponentials: if $P$ is a parabolic subgroup of a reductive group $G$,
and $p$ exceeds the nilpotence class $n(P)$ of the unipotent radical
$\V$ of $P$, then there is a $P$-equivariant isomorphism $\vlie \to
\V$ of algebraic groups, where $\vlie = \Lie(\V)$ is regarded as an
algebraic group via the Hausdorff formula.  As a consequence, we
observe in \ref{sub:intrinsic} that one gets an intrinsically defined
morphism of $P$-varieties $\NN(d,\vlie) \to \A(d,\V)$ which we prove
is injective. We have not so far been able to decide whether this
morphism should be an isomorphism of varieties, or even surjective.

\subsection{Some notations and conventions}

If $\Lambda$ is any commutative ring, and $V$ a finitely generated
$\Lambda$-module, we denote by $\Aut_\Lambda(V)$ the linear
automorphisms of $V$.  We denote by $\GL(V)$ the \emph{affine group
  scheme of finite type} with $\GL(V)(\Lambda') = \Aut_{\Lambda'}(V
\tensor_\Lambda \Lambda')$ for each commutative $\Lambda$-algebra
$\Lambda'$.

If $\Lambda = \EE$ is a field, we mostly prefer to identify affine
group schemes of finite type over $\EE$ which are absolutely reduced
with the corresponding linear algebraic groups over $\EE$.  Thus a
finite dimensional $\EE$-vector space $V$ determines a linear
algebraic group $\GL(V)$ over $\EE$; this is an $\EE$-form of $\GL(V
\tensor_\EE \overline{\EE})$, where $\overline{\EE}$ denotes an
algebraic closure of $\EE$.

\subsection{}
I would like to thank {\AA}rhus University for its hospitality during
the academic year 2000/1.  Also, thanks to Jens Jantzen for a number
of helpful comments on this manuscript.

\section{Nilpotent endomorphisms, groups, and Lie algebras}
\label{sec:nilpotent}

If $M$ is an Abelian group and $\phi$ is an nilpotent endomorphism of $M$,
the nilpotence degree of $\phi$ is the least positive integer $e$
such that $\phi^e = 0$.

A group $M$, respectively a Lie algebra $M$, is nilpotent provided
that its descending central series $M = C^0M \supseteq C^1M
\supseteq \cdots$ terminates in 1, respectively 0, after finitely many
steps [recall that for $i \ge 1$, we have $C^iM = (M,C^{i-1}M)$
for a group $M$, respectively $C^iM = [M,C^{i-1}M]$ for a Lie
algebra $M$]. If $M$ is nilpotent, its nilpotence class is the least
$e$ for which $C^eM$ is trivial.

Let $k$ be an algebraically closed field, and let $L$ be the Lie
algebra of a linear algebraic $k$-group; there is then a well-defined
notion of a nilpotent element of $L$. Suppose now that the
characteristic of $k$, say $p$, is positive; then $L$ is a $p$-Lie
algebra.  Evidently $X \in L$ is nilpotent if and only if $X^{[p^e]} =
0$ for some $e$ [the map $X \mapsto X^{[p^e]}$ is the $e$-th iteration
of the $p$-power map on $L$]. The $p$-nilpotence degree of a nilpotent
$X \in L$ is the minimal $e$ for which $X^{[p^e]} = 0$.  The element
$X$ is said to be $p$-nilpotent if its $p$-nilpotence degree is 1
(i.e. if $X^{[p]} = 0$).

We have (see \cite[\S3.1]{Bor1}) for all $X_1,X_2,\dots,X_r \in L$:
\begin{equation}
  \label{eq:Jacobson}
  \left(\sum_{i=1}^r X_i\right)^{[p]} \congruent \sum_{i=1}^r X_i^{[p]} \pmod{C^pL}
\end{equation}

If moreover $L$ is a nilpotent Lie algebra, then each element $X \in
L$ is nilpotent.  Since $L$ is a finite dimensional $p$-Lie algebra,
the $p$-nilpotence degree of any $X$ in $L$ is bounded; call the
$p$-exponent of $L$ the maximum of the $p$-nilpotence degree of its
elements.

\subsection{}
We have the following general result bounding ``exponent'' in terms of
``nilpotence class''.
\label{sub:general_nilp_lemma}
\begin{lem}
  \begin{enumerate}
  \item[(a)] Let $L$ be a nilpotent Lie algebra with class $e$. Assume
    for each $i \ge 0$ that $C^iL$ has a $k$-basis of $p$-nilpotent
    elements. Then $X^{[p^e]} = 0$ for all $X \in L$ (i.e.  the
    $p$-exponent of $L$ is $\le e$).
  \item[(b)] Let $G$ be a nilpotent group with class $e$. Assume for
    each $i \ge 0$ that $C^iG$ is generated by elements of order $p$.
    Then $x^{p^e} =1$ for all $x \in G$ (i.e. the exponent of $G$ is
    $\le p^e$).
  \end{enumerate}
\end{lem}

\begin{proof}
  Note that for $X=G$ or $X=L$ we have $C^iC^jX \subset C^{ij}X$ for
  all $i,j \ge 1$. To prove (a), let $X \in L$ and write $X =
  \sum_{i=1}^r X_i$ where $X_i^{[p]} = 0$ for all $i$.  Then $X^{[p]}
  \in C^pL$ by \eqref{eq:Jacobson}.  We have $C^{p^{e-1}}C^pL = 0$ by
  assumption, so by induction $X^{[p^e]} = (X^{[p]})^{[p^{e-1}]} = 0$.
  The proof of (b) is essentially the same.
\end{proof}

\section{Abelian unipotent groups}
\label{sec:Abelian-unipotent}

In this section, we recall some basic known facts about Abelian
unipotent groups.

\subsection{Witt vector groups.}
\label{sub:witt}
Let $p > 0$ be a prime number, let $n \ge 1$ be an integer, and let
$\Witt_{n,\Zp}$ denote the group scheme over $\Zp$ of the ``Witt
vectors of length $n$'' for the prime $p$; see
\cite[II.\S6]{serre79:_local_field} and
\cite[V.\S16,VII.\S7]{SerreAG&CF}.  We write $\Witt_n = \Witt_{n,k}$
for the corresponding group over $k$.

\begin{example}
  Let $F(X,Y) =  \dfrac{X^p + Y^p - (X + Y)^p}{p} \in \Z[X,Y]$. When $n
  = 2$, the operation in $\Witt_{2,\Zp} = \Aff^2_{\Zp}$ (here written
  additively) is defined by the rule
  \begin{equation*}
    \vec{t} + \vec{s} = (t_0 + s_0,
    F(t_0,s_0) + t_1 + s_1).
  \end{equation*}
  More precisely, the co-multiplication for $\Zp[\Witt_n] = \Zp[T_0,T_1]$
  is given by 
  \begin{equation*}
    \Delta(T_0) = T_0 \tensor 1 + 1 \tensor T_0
  \end{equation*}
  and
  \begin{equation*}
    \Delta(T_1) = T_1 \tensor 1 + 1 \tensor T_1 + 
    F(T_0 \tensor 1,1 \tensor T_0)
  \end{equation*}
\end{example}

\begin{stmt}
  \label{stmt:witt-free-coord-ring}
  The underlying scheme of \ $\Witt_{n,\Zp}$ is isomorphic with the
  affine space ${\Aff^n}_{\Zp}$, hence the structure algebra
  $\Zp[\Witt_n]$ is free over $\Zp$.
\end{stmt}

\begin{stmt}
  \label{stmt:witt-char-0}
  There is an isomorphism
  of \ $\Q$-group schemes
  \begin{equation*}
    \varphi:\Witt_{n,\Q} \xrightarrow{\iso} \G_{a,\Q} \times \cdots \times
    \G_{a,\Q} \quad \text{($n$ factors).}
\end{equation*}
\end{stmt}

\begin{proof}
  For $m \ge 1$, let 
  \begin{equation*}
    w_m = X_0^{p^m} + pX_1^{p^{m-1}} + \cdots + p^mX_m
    \in \Z[X_0,X_1,\dots].
  \end{equation*}

  We may define a map
  \begin{equation*}
    \varphi:\Witt_{n,\Q} \to \G_{a,\Q} \times \cdots \times  \G_{a,\Q}
  \end{equation*}
  by assigning, for each $\Q$-algebra $\Lambda$ and each $\vec{t} \in
  \Witt_{n,\Q}(\Lambda)$, the value
  \begin{equation*}
    \varphi(\vec{t}) =
    (w_0(\vec{t}),w_1(\vec{t}),\dots,w_{n-1}(\vec{t})).
  \end{equation*}    
  Since $p$ is invertible in $\Q$, it follows from \cite[Theorem
  II.6.7]{serre79:_local_field} that $\varphi$ is an isomorphism of
  $\Q$-group schemes. [Note that the assertion is valid over any field
  $\FF$ provided only that the characteristic of the field $\FF$ is
  different from $p$.]
\end{proof}

\subsection{The Artin-Hasse exponential series}
\label{sub:Artin-Hasse-series}
Now let $F(t) \in \Q \pow{t}$ be the power series
\begin{equation*}
  F(t) = \exp(-(t + t^p/p + t^{p^2}/p^2 + \cdots)).
\end{equation*}
If $\mu$ denotes the M\"obius function, one easily checks the identity
of formal series
\begin{equation*}
  F(t) = \prod_{(m,p)=1, m \ge 1} (1-t^m)^{\mu(m)/m},
\end{equation*}
by taking logarithms and using the fact that
\begin{equation*}
  \sum_{d|m} \mu(d) = 0
\end{equation*} 
if $m \not = 1$. It then follows that \emph{the coefficients of $F(t)$
  are integers at $p$}; i.e. $F(t) \in \Z_{(p)}\pow{t}$.

\subsection{}
\label{sub:Artin-Hasse-linear}

If $\LL$ is a $\Z_{(p)}$-lattice, and $X \in \End_{\Z_{(p)}}(\LL)$ is
a nilpotent endomorphism such that $X^{p^n} = 0$, then there is a
homomorphism of $\Z_{(p)}$-group schemes
\begin{equation*}
  E_X:\Witt_{n,\Z_{(p)}} \to \GL(\LL)
\end{equation*}
given for each $\Z_{(p)}$-algebra $\Lambda$ by $E_X(\vec{t}) =
F(t_0X)F(t_1X^p)\cdots F(t_{n-1}X^{p^{n-1}})$ for $\vec{t} \in
\Witt_n(\Lambda)$; see \cite[V\S16]{SerreAG&CF}.

Let $V_\Q = \LL \tensor_{\Zp} \Q$.
There are maps
\begin{equation*}
  E_{X}:\Witt_{n,\Q} \to \GL(V_\Q) \quad \text{and} \quad
  E_{\overline{X}}:\Witt_{n,\F_p} \to \GL(\LL/p\LL).
\end{equation*}
obtained by base change.

\begin{lem}
  \begin{enumerate}
  \item Over $\Q$, $E_{X}$ factors as
    \begin{equation*}
    \xymatrix{
      \Witt_{n,\Q} \ar[r]^{E_X} \ar[dr]_\varphi  & \GL(V_\Q) \\
      & \G_{a,\Q} \times \cdots \times \G_{a,\Q} 
      \ar[u]_{\vec{s} \mapsto \exp(\sum_{j=0}^{n-1} p^{-j}s_jX^{p^j})} \\
      }
  \end{equation*}
  where $\varphi$ is the isomorphism of \ref{sub:witt}\ref{stmt:witt-char-0}.
\item The endomorphism $\overline{X}$ of $\LL/p\LL$ is in the image of
  the Lie algebra homomorphism
  \begin{equation*}
    dE_{\overline{X}}:\wlie_{n,\F_p} \to \lie{gl}(\LL/p\LL).
  \end{equation*}
\end{enumerate}
\end{lem}

\begin{proof}
  For each $\Q$-algebra $\Lambda$ and each $\vec{t} \in
  \Witt_n(\Lambda)$ one uses induction on $n$ and the definition of
  the $w_j$ to verify that
  \begin{equation*}
    \sum_{j=0}^{n-1} p^{-j}w_j(\vec{t})X^{p^j}  = 
    \sum_{m=0}^{n-1}\sum_{l=0}^{n-1-m}  p^{-l} (t_m X^{p^m})^{p^l}.
  \end{equation*}
  It follows that
  \begin{equation*}
    \exp\left(\sum_{j=0}^{n-1} p^{-j}w_j(\vec{t})X^{p^j} \right) = 
    \prod_{m=0}^{n-1} F(t_mX^{p^m}) = E_X(\vec{t}),
  \end{equation*}
  whence (1).
  
  For (2), let $T_0,\dots,T_{n-1}$ denote the coordinate functions on
  $\Witt_{n,\F_p}$ with $T_i(\vec{t}) = t_i$; thus $A=\F_p[\Witt_n]$ is a
  polynomial ring in the $T_i$.  The tangent space to $\Witt_n$ at $0$
  contains the ``point-derivation'' $D:A \to \F_p$ given by $f \mapsto
  \dfrac{\partial f}{\partial T_0}\mid_{\vec{t} = 0}$, and it is clear
  that $dE_{\overline{X}}(D) = \overline{X}$.
\end{proof}

\subsection{The Lie algebra of the Witt vectors}
\label{sub:Witt-Lie}

Let $V$ be a $k$-vector space of dimension $p^{n-1} + 1$, and let $X
\in \End_k(V)$ be a nilpotent ``Jordan block'' of size $p^{n-1}+1$.
Thus $x^{p^{n-1}} \not = 0$ while $x^{p^n} = 0$.  The smallest $p$-Lie
subalgebra of $\lie{gl}(V)$ containing $x$ is then the Abelian Lie
algebra $\lie{a} = \sum_{i=0}^{n-1} kX^{p^i} = \sum_{i=0}^{n-1} k
X^{[p^i]}.$

Let $E_X:\Witt_n \to \GL(V)$ be the homomorphism determined by $X$
as in \ref{sub:Artin-Hasse-linear}.

\begin{prop}
  \label{stmt:Witt-Lie-algebra}
  If $k$ is a field of characteristic $p$, then $\wlie_n =
  \Lie(\Witt_n)$ is an Abelian Lie algebra with a $k$-basis
  $Z_0,Z_1,\dots,Z_{n-1}$ such that $Z_i = Z_0^{[p^i]}$ for $0 \le i
  \le n-1$.
\end{prop}
\begin{proof}
  Let $\lie{b}$ be the image of $dE_X$; thus $\lie{b}$ is a $p$-Lie
  subalgebra of $\lie{gl}(V)$.  According to Lemma
  \ref{sub:Artin-Hasse-linear}, $\lie{b}$ contains $X$. It follows
  that $\lie{a} \subset \lie{b}$.  On the other hand, we have
  \begin{equation*}
    n = \dim_k \lie{a} \le \dim_k \lie{b} \le \dim \Witt_{n,k} = n.
  \end{equation*}
  Thus $\lie{a} = \lie{b} \iso \wlie_n$, and the proposition follows.
\end{proof}

\begin{example}
  Say $n=2$. Then $k[\Witt_2] = k[T_0,T_1]$. One can show that $X_0 =
  \dfrac{\partial}{\partial T_0} + T_0^{p-1}\dfrac{\partial}{\partial
    T_1}$ and $X_1 = \dfrac{\partial}{\partial T_1}$ are
  $\Witt_2$-invariant derivations of $k[\Witt_2]$, and that these
  derivations span $\wlie_2$.  A simple computation yields $X_0^{[p]}
  = X_1$.
\end{example}

\subsection{Exponents} 

\begin{prop}
  \label{stmt:Witt-exponent}
  \begin{enumerate}
  \item The $p$-exponent of $\wlie_n = \Lie(\Witt_n)$ is $n$.
  \item The exponent of $\Witt_n$ is $p^n$; moreover, a Witt vector
    $\vec{t} \in \Witt_n(k)$ has order $p^n$ if and only if $t_0 \not
    = 0$.
  \end{enumerate}
\end{prop}
\begin{proof}
  Part (1) follows immediately from the description of $\wlie_n$ given
  by Proposition \ref{sub:Witt-Lie}.  For (2), recall \cite[Theorem
  II.6.8]{serre79:_local_field} that the ring of (infinite) Witt
  vectors $\Witt(k)$ is a strict $p$-ring (see \emph{loc. cit.} II.5
  for the definition) and that $\Witt_m(k) \iso \Witt(k)/p^m\Witt(k)$
  for all $m \ge 1$. (2) now follows at once.
\end{proof}

\subsection{Connected Abelian unipotent groups}
\label{sub:connect-abelian}

Recall that two connected Abelian algebraic groups $G$ and $H$ are
said to be \emph{isogenous} if there is a surjection $G \to H$ whose
kernel is finite.

\begin{lem}
  If $G$ and $H$ are connected Abelian algebraic groups which are
  isogenous, then the exponent of $G$ is equal to the exponent of $H$.
\end{lem}
\begin{proof}
  The lemma is clear if either $G$ or $H$ has infinite exponent, so
  assume otherwise.  Suppose that $\phi:G \to H$ is a surjection with
  finite kernel. Let $m$ be the exponent of $H$.  Since $G$ is
  Abelian, the map $x \mapsto x^m$ defines a group homomorphism $G \to
  \ker \phi$; since $G$ is connected and $\ker \phi$ is finite, this
  homomorphism must be trivial. It follows that $x^m = 1$ for all $x
  \in G$, and this shows that the exponent of $G$ is $\le$ that of
  $H$. The inequality $\ge$ is immediate since $\phi$ is surjective.
\end{proof}

\begin{prop}Let $U$ be a connected Abelian unipotent group over $k$. Then
  \begin{enumerate}
  \item $U$ is isogenous to a product of Witt groups $\prod_{i=1}^d
    \Witt_{n_i,k}$; moreover, the integers $n_i$ are uniquely
    determined (up to order) by $U$.
  \end{enumerate}
  Let $n = \max_i(n_i)$ where the $n_i$ are as in 1.
  \begin{enumerate}
  \item  The exponent of the group $U$ is $p^n$.
  \item  The $p$-exponent of $\ulie = \Lie(U)$ is $\le n$.
  \end{enumerate}
\end{prop}

\begin{proof}
  The first assertion is \cite[VII\S 2 Theorem 1]{SerreAG&CF}. The
  second assertion follows immediately from the lemma.
  
  For the last assertion, it is proved in \cite[VII\S 2 Theorem
  2]{SerreAG&CF} that the group $U$ is a \emph{subgroup} of a product
  of Witt groups. A careful look at the proof in \emph{loc. cit.}
  shows that the exponent of $U$ and this product may be chosen to
  coincide.  Thus, $\ulie$ is a subalgebra of $\wlie$, a product of
  Lie algebras $\wlie_{n_i}$ with $\max(p^{n_i})$ equal to the
  exponent of $U$; (3) now follows since the $p$-exponent of
  the Lie subalgebra $\ulie$ can't exceed that of $\wlie$.
\end{proof}

\begin{rem}
  The $p$-exponent of $\Lie(U)$ may indeed be strictly smaller than
  $\log_p$ of the exponent of $U$.  Let $\mathrc{V}_2$ be the
  algebraic $k$-group which is isomorphic as a variety to $\Aff^2_k$,
  with the group operation in $\mathrc{V}_2$ determined by
  \begin{equation*}
    \vec{t} + \vec{s} = (t_0 + s_0,
    F(t_0,s_0)^p + t_1 + s_1).
  \end{equation*}
  Then the map $\varphi:\Witt_2 \to \mathrc{V}_2$ given by $\vec{t}
  \mapsto (t_0^p, t_1)$ is a (purely inseparable) isogeny. The
  exponent of $\mathrc{V}_2$ is $p^2$, but every element $x \in
  \vlie_2 = \Lie(\mathrc{V}_2)$ satisfies $x^{[p]} = 0$. Indeed, one
  can check that $\dfrac{\partial}{\partial T_0}$ and
  $\dfrac{\partial}{\partial T_1}$ are $\mathrc{V}_2$-invariant
  derivations of $k[\mathrc{V}_2]$, and that they span $\vlie_2$ over
  $k$.
\end{rem}

\section{Reductive groups}
\label{sec:reductive}

\subsection{Generalities}

Let $G$ be a connected reductive group over the field $k$ which is
defined and split over the prime field $\F_p$.  We fix a maximal torus
$T$ contained in a Borel subgroup $B$ of $G$.  Let $X = X^*(T)$ be the
group of characters of $T$, and $Y= X_*(T)$ be the group of
co-characters.  The adjoint action of $G$ on its Lie algebra $\glie$
is diagonalizable for $T$; the non-zero weights of this action form a
root system $R \subset X$, and the choice of Borel subgroup determines
a system of positive roots $R^+$ and a system of simple roots $S$.
Write $\langle ?,? \rangle$ for the canonical pairing $X \times Y \to
\Z$.



For each root $\alpha \in R^+$, there is a root homomorphism
$\phi_\alpha:\G_a \to U$; the subgroup $U$ is equal to the direct
product (in any fixed order) of the images of the root homomorphisms
$\phi_\alpha$ with $\alpha > 0$.

For each $\alpha \in R^+$ the derivative of $\phi_\alpha$ yields an
element $e_\alpha \in \ulie = \Lie(U)$; the $e_\alpha$ form a
basis for $\ulie.$

\subsection{Good primes}

We will usually assume that $p$ is a good prime for $G$. If the root
system of $G$ is indecomposable, let $\beta$ be the short root of
maximal height. In that case, the prime $p$ is good for $G$ provided
that if $\beta^\vee = \sum_{\alpha \in S} a_\alpha \alpha^\vee$, then
all $a_\alpha$ are prime to $p$. For indecomposable root systems, $p$
is bad (=not good) just in case one of the following holds: $p=2$ and
$R$ is not of type $A_r$; $p=3$ and $R$ is of type $G_2$, $F_4$ or
$E_r$; or $p=5$ and $R$ is of type $E_8$. In general $p$ is good for
$G$ if it is good for each indecomposable component of the root system
$R$.

\subsection{Parabolic subgroups}
\label{sub:parabolic}

Let $P$ be a parabolic subgroup of $G$ containing the Borel subgroup
$B$, and let $\plie$ be the Lie algebra of $P$. Put $I = \{\alpha \in
S \mid \plie_{-\alpha} \not = 0\}$. The parabolic subgroup $P$ is then
\begin{equation*}
  P = \langle B,\imag\phi_{-\alpha} \mid \alpha \in I \rangle.
\end{equation*}
The group $P$ has a Levi decomposition $P= \Levi\V$ where $\Levi$ is a
reductive group and $\V$ is the unipotent radical of $P$.  The derived
group of the Levi factor $\Levi$ is a semisimple group whose root
system $R_P$ is generated by the roots in $I$.  Denote by $\vlie =
\Lie(\V)$ the nilradical of $\plie$. The group $\V$ is the product of
the images of the root homomorphism $\phi_\alpha$ with $\alpha \in R^+
\setminus R_P$.

There is (see e.g. \cite[Ch.  9]{springer98:_linear_algeb_group}) an
isogeny
\begin{equation*}
  \hat G =\prod_i G_i \times T \to G
\end{equation*}
where each $G_i$ is semisimple with indecomposable root system, and
$T$ is a torus.  Let $\hat P$ denote the parabolic subgroup of $\hat
G$ determined by $I$, and let $\hat \V$ denote its unipotent radical.
\begin{lem}
  The above isogeny restricts to an isomorphism $\hat \V \iso \V$;
  moreover, we have $\hat \V \iso \prod_i \hat \V_i$ where $\hat \V_i
  = \V \cap G_i$.
\end{lem}

\begin{proof}
  This follows from \cite[Prop. 22.4]{Bor1}.
\end{proof}

\subsection{}
\label{sub:nilp-class}

Associated with the parabolic subgroup $P$, we may define a
homomorphism $f:\Z R \to \Z$ given by
\begin{equation*}
  f(\alpha) = \left \{ 
    \begin{matrix}
      0 & \text{if $\alpha \in S$ and $-\alpha \in R_P$} \\
      2 & \text{if $\alpha \in S$ and $-\alpha \not \in R_P$}
    \end{matrix}
\right .
\end{equation*}
Such a homomorphism induces a grading of the Lie algebra $\glie =
\bigoplus_{i \in \Z} \glie(i)$ by setting $\glie(i) =
\bigoplus_{f(\alpha)=i} \glie_\alpha$. We have evidently $\Lie(P) =
\plie = \bigoplus_{i \ge 0}\glie(i)$ and $\Lie(\V) = \vlie =
\sum_{i>0}\glie(i)$. We have by construction that 
\begin{equation*}
  \glie(i) \not = 0 \implies 
      i \congruent 0 \pmod 2 .
\end{equation*}

When $R$ is indecomposable, let $\tilde \alpha \in R^+$ be the long
root of maximal height, and let $n(P) = \frac{1}{2}f(\tilde \alpha) +
1$. If we write $\tilde \alpha$ as a $\Z$-linear combination of the
simple roots $S$, then $n(P) -1 $ is just the sum of the coefficients
in this expression of the roots in $S \setminus I$.

Note that 
\begin{equation*}
  \glie(i) \not = 0 \implies -f(\tilde \alpha) \le i \le f(\tilde \alpha)
\end{equation*}
and that $\Lie(T) \subset \glie(0) \not = 0$ and $e_{\tilde \alpha}
\in \glie(f(\tilde \alpha)) \not= 0$.

When $R$ is no longer indecomposable, let $S'$ be the simple roots for
an indecomposable component $R'$ of $R$, and let $\tilde \alpha'$ be
the highest long root in $R'$. Put $n(P,S') = \dfrac{1}{2}f(\tilde
\alpha') +1$, and let $n(P)$ be the supremum of the $n(P,S')$.

\begin{prop} Suppose that $p$ is a good prime for $G$, let $P$ be
  a distinguished parabolic subgroup, and let $m \ge 1$ be minimal
  with $p^m \ge n(P)$.
  \begin{enumerate}
  \item[(a)] The nilpotence class of the Lie algebra $\vlie$ is $n(P)$.
  \item[(b)] The $p$-exponent of $\vlie$ is $\le m$.
  \item[(c)] The nilpotence class of the group $\V$ is $n(P)$.
  \item[(d)] The exponent of $\V$ is $\le p^m$.
  \end{enumerate}
\end{prop}

\begin{proof}
  We first prove (a) and (c).  By Lemma \ref{sub:parabolic}, we are
  reduced to the case where $R$ is indecomposable.
  
  Since $p$ is good, \cite[Prop.  4.7]{Borel-Tits} shows that
  $C^{j-1}\V = \prod_{f(\alpha) \ge 2j} \imag\phi_\alpha$ for $j
  \ge 1$.  Essentially the same arguments show that $C^{j-1}\vlie =
  \sum_{f(\alpha) \ge 2j} ke_\alpha$. Since every root $\alpha$
  satisfies $f(\tilde \alpha) \ge f(\alpha)$; (a) and (c) now follow
  at once.
  
  Note that we have showed for all $i \ge 0$ that $C^i\vlie$ has a
  basis of $p$-nilpotent vectors (the root vectors) and that $C^i\V$
  is generated by elements of order $p$ (the images of root
  homomorphisms). By Lemma \ref{sec:nilpotent}, (b) now follows from
  (a), and (d) follows from (c).
\end{proof}

\begin{cor}
  If $p \ge h$, every nilpotent element $Y \in \glie$ satisfies
  $Y^{[p]}=0$ and every unipotent element $1 \not = u \in G$ has
  order $p$.
\end{cor}

\begin{proof}
  Let $Y$ be a \emph{regular nilpotent} element in $\ulie = \Lie(U)$;
  thus $Y$ is a representative for the dense $B$-orbit on $\ulie$ [see
  the discussion of Richardson's dense orbit theorem below in section
  \ref{sub:bala-carter}]. Since $n(B) = h-1$, the proposition shows
  that $Y^{[p]} = 0$. Since the regular nilpotent elements form a
  single dense orbit in the nilpotent variety, we get $Y^{[p]} = 0$
  for every nilpotent $Y$.

  The assertion for unipotent elements follows in the same way.
\end{proof}

\begin{rems}
  \begin{enumerate}
  \item Suppose that $G$ is semisimple; thus $S$ is a $\Q$ basis for
    $X_\Q$. Then $f$ determines, by extension of scalars, a unique
    $\Q$-linear map $f_\Q:X_\Q \to \Q.$ For some $q \in \Z$, the
    homomorphism $qf_\Q$ maps $X$ to $\Z$, so that $qf = \phi \in Y$,
    the group of cocharacters of the maximal torus $T$. [In fact, this
    is true even when $G$ is only assumed to be reductive, rather than
    semisimple.]  For any such integer $q$, we have $\glie(i) = \{Y
    \in \glie \mid \Ad(\phi(t))Y = t^{qi}Y\}$, which makes it clear
    that the grading of $\glie$ is a grading as a $p$-Lie algebra. In
    particular, if $Y \in \glie(i)$, we have $Y^{[p]} \in \glie(pi)$.
  \item The preceding remark permits an alternate proof of (b) of the
    proposition; indeed, for a homogeneous $Y \in \vlie(i)$ [so
    $i>0$], we have $Y^{[p^m]} \in \glie(p^mi) = 0$.
  \item The corollary is of course well known; I didn't find a
    suitable reference, however. Jens Jantzen has pointed out to me a
    somewhat more elementary argument that $X^{[p]} = 0$ for $X \in
    \glie$ nilpotent when $p \ge h$.  We may assume $X$ to be in
    $\ulie$; thus we may write $X = \sum_{\alpha \in R^+} a_\alpha
    e_\alpha$ with scalars $a_\alpha \in k$.  By Jacobson's formula
    for the $p$-th power of a sum, $X^{[p]}$ is $\sum_{\alpha \in R^+}
    a_\alpha^p e_\alpha^{[p]} + L$ where $L$ is a linear combination
    of commutators of length $p$.  Now, all summands of $L$ are weight
    vectors of a weight that has height $\ge p$. But the maximal
    height of a root is $h - 1 < p$.
  \end{enumerate}
\end{rems}

\subsection{The Bala-Carter parameterization of nilpotent elements}
\label{sub:bala-carter}

Let $G$ be connected and reductive in good characteristic $p$, and let
$P = \Levi\V$ be a parabolic subgroup.  The adjoint action $P$ on
$\glie$ leaves $\vlie$ invariant. A theorem of Richardson
\cite[Theorem 5.3]{hum-conjugacy} guarantees that $P$ has a unique
open orbit on $\vlie$, and a unique open orbit on $\V$; these are the
Richardson orbits of $P$, and representatives for these orbits are
called Richardson elements.

A nilpotent element $X \in \glie$ is \emph{distinguished} if the
connected center of $G$ is a maximal torus of $C_G(x)$. [If $G$ is
semisimple, this means that any semisimple element of $C_G(x)$ is
central].  

On the other hand, the parabolic subgroup $P$ is called
\emph{distinguished} if
\begin{equation*}
  \dim \glie(0) - \dim Z(G) = \dim \glie(2).
\end{equation*}
[Note that this differs from the definition in \cite[p.167]{Carter1},
but that Corollary 5.8.3 of \emph{loc. cit.} shows that it is
equivalent in case $p$ is good.]

The following relates these two notions of distinguished:
\begin{prop}  \cite[Prop. 5.8.7]{Carter1} If $P$ is
  distinguished, then a Richardson element in $\vlie$ is a
  distinguished nilpotent element. Moreover, the Richardson orbit on
  $\vlie$ meets $\glie(2)$ in an open $\Levi$-orbit.
\end{prop}

The full Bala-Carter theorem is as follows:
\begin{prop}\emph{Bala,Carter \cite{Bala-CarterI,Bala-CarterII},
    Pommerening \cite{PommereningI,PommereningII}} There is a
  bijection between the $G$-orbits of nilpotent elements in $\glie$
  and the conjugacy classes of pairs $(L,Q)$ where $L$ is a Levi
  subgroup of a parabolic subgroup of $G$, and $Q$ is a distinguished
  parabolic subgroup of $L$. The nilpotent orbit determined by $(L,Q)$
  is the one meeting the nilradical of $\Lie(Q)$ in its Richardson
  orbit.
\end{prop}

\section{The $p$-exponent of $\vlie$}
\label{sec:Lie-alg}

Throughout this section and the next, $G$ is a reductive group, $P$ is
a distinguished parabolic subgroup with Levi decomposition $P = \Levi
\V$. The characteristic $p$ is assumed to be good for $G$.

\subsection{}
\label{sub:order-generalities}
Let $A$ be a discrete valuation ring, with residue field of
characteristic $p$ and quotient field $\FF$.  We assume chosen some
fixed embedding of the residue field of $A$ in our algebraically
closed field $k$.

If $\LL$ is an $A$-lattice and $\psi$ is a nilpotent $A$-endomorphism,
one might hope to relate the Jordan block structure of $\psi_k$ and
$\psi_\FF$ [If $\psi \in \End_A(\LL)$, the corresponding endomorphisms
of $L_\FF = \LL \tensor_A \FF$ and $L_k$ will be denoted $\psi_\FF$
and $\psi_k$]. Since the dimension of the kernel of a liner
transformation is equal to the number of its Jordan blocks, one must
require that $\dim_k \ker \psi_k = \dim_\FF \ker \psi_\FF$. However,
even with that condition, the partitions can be different; indeed, let
$\pi$ be a prime element of $A$ and consider the endomorphism $\psi$
of the lattice $A^4 = \bigoplus_{i=1}^4 A e_i$ determined by the rules
$\psi(e_1) = 0$, $\psi(e_2) = \psi(e_3) = e_1$, $\psi(e_4) =
(\pi-1)e_2 + e_3$.  Then $\psi_\FF$ has partition $(3,1)$ while
$\psi_k$ has partition $(2,2)$.

On the other hand, one has the following straightforward result. Let
$\LL$ be an $A$-lattice which is $2\Z$-graded; say
\begin{equation*}
  \LL = \bigoplus_{i=0}^d \LL_{2i}, 
  \quad \LL_0 \not = 0, \quad \LL_{2d} \not = 0.
\end{equation*}\
Let $\LL^+ = \bigoplus_{i>0} \LL_{2i}$, $L_\FF^+$ and $L_k^+$ be in
each case the sum of the homogeneous components of positive degree.

\begin{prop}
  Suppose $\psi \in \End_A(\LL)$ is an endomorphism of degree 2 (i.e.
  $\psi(\LL_{i}) \subseteq \LL_{i+2}$ for all $i$), and assume
  $\psi_\FF:L_\FF \to L^+_\FF$ is surjective, so that the nilpotence
  degree of $\psi$ is $d+1$.  If $\dim_k \ker \psi_k = \dim_\FF \ker
  \psi_\FF$, then $\psi:\LL \to \LL^+$ is surjective.  In particular,
  $\psi_k^d \not = 0$, hence the nilpotence degree of $\psi_k$ is also
  $d+1$.
\end{prop}

\begin{proof}
  It suffices to show that $\psi_k:L_k \to L^+_k$ is surjective. Since
  $\dim_\FF L_\FF = \dim_k L_k$ and $\dim_\FF L_\FF^+ = \dim_k L_k^+$,
  that follows immediately from the assumption on kernel dimensions.
\end{proof}

\subsection{}
\label{sub:spaltenstein}

Let $G_{\Z}$ be a split reductive group scheme over $\Z$ which gives
rise to $G$ upon base change.  There is a general notion of the Lie
algebra $\glie_{\Z}$ of the affine group scheme $G_{\Z}$; see
\cite[I.7.7]{JRAG}.  In the case of our split reductive group,
$\glie_{\Z}$ may be described explicitly; see \cite[II.1.11]{JRAG}.  For
any commutative ring $A$, let $\glie_{A} = \glie_{\Z} \tensor_\Z A$.  The
explicit description of $\glie_{\Z}$ implies that it is a $\Z$-lattice in
the split reductive $\Q$-Lie algebra $\glie_{\Q}$, and that $\Lie(G) =
\glie = \glie_{k}$.  If $X_\Z \in \glie_{\Z}$, denote by $X_A$ the element
$X_\Z \tensor 1 \in \glie_{A}$.

Let $P$ be a distinguished parabolic subgroup $P$, and let the map
$f:\Z R \to \Z$ be as in \ref{sub:nilp-class}; then $f$ induces also a
grading of $\glie_{\Z}$ (this again relies on the explicit description
of $\glie_{\Z}$ mentioned above).

Fix an algebraic closure $\overline{\Q}$ of $\Q$.
\begin{lem}
  There is a finite subextension $\Q \subset \FF \subset
  \overline{\Q}$ and a valuation ring $A \subset \FF$, whose residue
  field we embed in $k$, such that the following holds: for some $X_A
  \in \glie_A(2)$, the element $X_k$ is Richardson  in
  $\vlie$, and $X_\FF$ is Richardson in
  $\vlie_{\overline{\Q}}$.
\end{lem}

[This lemma is implicit in \cite{spalt-transverse}; we include a proof
for the convenience of the reader.]

\begin{proof}
  The Richardson orbit on $\vlie_{\bar \Q}$ meets $\glie_{\bar \Q}(2)$
  in an open set, so we may find a regular function $f$ on
  $\glie_{\bar \Q}(2)$ such that $f(Y) \not = 0$ implies $Y$ is a
  Richardson element.  Using the lattice $\glie_\Z(2)$ in
  $\glie_{\bar \Q}(2)$, we obtain coordinate functions [dual to the
  root-vector basis] on $\glie_{\bar \Q}(2)$; let $\FF \subset \bar
  \Q$ by a finite extension of $\Q$ containing the coefficients of $f$
  with respect to these coordinate functions.  Take for $A$ the
  localization of the ring of integers of $\FF$ at some prime lying
  over $(p)$, and fix some embedding of the residue field of $A$ in
  $k$.  If $\pi$ denotes a prime element of $A$, we may multiply $f$
  be a suitable power of $\pi$ and assume that all the coefficients of
  $f$ are in $A$, and that not all are 0 modulo $\pi$. Let $\hat f \in
  k[\glie(2)]$ be the function obtained by reducing $f$ modulo $\pi$.
  Then the distinguished open set determined by $\hat f$ is non-empty
  and so must meet the set of Richardson elements in $\glie(2)$.
  After possibly enlarging $\FF$ and $A$, we may suppose that there is
  a Richardson element $X \in \glie(2)$ with $\hat f(X) \not = 0$, and
  such that the coefficients of $X$ (in the root-vector basis) lie in
  the residue field of $A$.  It is then clear that any lift $X_A$ of
  $X$ to $\glie_A(2)$ has the desired property.
\end{proof}

The main result obtained by Spaltenstein in \cite{spalt-transverse}
implies the following:
\begin{prop}
  Assume that the root system of $G$ is indecomposable, and moreover
  that $p$ is a good prime if $R$ is not of type $A_r$, and that $G =
  \GL_{r+1}$ if $R=A_r$.  Choose a finite extension $\FF$ of $\Q$ with
  ring of integers $A$ as in the lemma; let $X_A \in \vlie_A$ be such
  that $X = X_k$ is a Richardson element of $\vlie$ and $X_\FF$ is a
  Richardson element of $\vlie_\FF$.  Then
  \begin{enumerate}
  \item $\clie_\glie(X) \subset \plie$, and
  \item $\dim_\FF \clie_{\glie_\FF}(X_\FF) = \dim_k \clie_\glie(X)$.
  \end{enumerate}
\end{prop}

\begin{rem}
  This result was also obtained by Premet in \cite{Premet--JacMor}.
\end{rem}

\subsection{}
\label{sub:sl2-lemma}

Let $\FF$ be a field of characteristic 0, and let $\sllie_2(\FF)$ be
the split simple Lie algebra over $\FF$ of $2 \times 2$ matrices with
trace 0. This Lie algebra has an $\FF$-basis $X,Y,H$, where $[H,X] =
2X$, $[H,Y] = -2Y$, and $[X,Y] = H$. The semisimple element $H$ acts
diagonally on any finite dimensional representation $(\rho,M)$, and
the weights of $H$ (=eigenvalues of $\rho(H)$) on $M$ are integers.
Write $M_i$ for the $i$-th weight space.

The following is an easy consequence of the classification of finite
dimensional representations for $\sllie_2(\FF)$:
\begin{lem}
  Let $(\rho,M)$ be a finite dimensional $\sllie_2(\FF)$-module such
  that all eigenvalues of $\rho(H)$ on $M$ are even. Then
  $\rho(X):\bigoplus_{i \ge 0} M_i \to \bigoplus_{i > 0} M_i$
  is surjective.
\end{lem}

\subsection{} We can now prove the statement of theorem 
\ref{sub:first_theorem} for the Lie algebra.

\label{sub:p-exponent}
\begin{theorem}
  \emph{($p$-exponent formula)} Let $m \ge 1$ be the minimal integer
  with $p^m \ge n(P)$. Then the $p$-exponent of $\vlie$ is $m$, and
  the $p$-nilpotence degree of a Richardson element of $\vlie$ is $m$.
\end{theorem}

\begin{proof}
  In view of Lemma \ref{sub:parabolic}, we may suppose that $G$
  satisfies the hypothesis of Proposition \ref{sub:spaltenstein}.
  
  With $m$ as in the statement of the theorem, Proposition
  \ref{sub:nilp-class} shows that the $p$-nilpotence degree of any
  element of $\vlie$ is $\le m$; to prove that equality holds, it
  suffices to exhibit a representation $(\rho,W)$ of $\plie$ as a
  $p$-Lie algebra in which some $\rho(X)^{n(P)-1}$ acts non-trivially.
  Take $(\rho,W) = (\ad,\plie)$; we show that $\ad(X)^{n(P)-1} \not =
  0$ for a suitable (and hence any) Richardson element.
  
  First, use the lemma to find a finite extension $\FF$ of $\Q$, a
  valuation ring $A \subset \FF$, and an element $X_A \in \glie_A(2)$
  such that $X=X_k \in \vlie$ is Richardson, and $X_\FF \in
  \vlie_{\overline{\Q}}$ is Richardson.
  
  The discussion in \ref{sub:nilp-class} shows that the lattice
  $\plie_A$ is $2\Z$ graded as in \ref{sub:order-generalities} with
  $d=n(P)-1$.  In the notation of Lemma \ref{sub:order-generalities},
  we have $L_\FF = \plie_\FF$ and $L_\FF^+ = \vlie_\FF$.  Note that
  $\ad(X_A):\plie_A \to \plie_A$ has degree 2.
  
  It follows from \cite[Prop. 5.8.8]{Carter1} (or more precisely, the
  proof of that Proposition) that there are elements $Y,H \in
  \glie_{\overline{\Q}}$ such that $H$ is semisimple, $\overline{\Q}Y
  + \overline{\Q}H + \overline{\Q}X_\FF$ is a subalgebra isomorphic to
  $\sllie_2$, and such that the grading of $\glie_{\overline{\Q}}$
  determined by $H$ is the same as that determined by the function $f$
  as in \ref{sub:nilp-class}. Thus for each $i \in \Z$ we have
  $\glie_{\overline{\Q}}(i) = \{Z \in \glie_{\overline{\Q}} \mid [H,Z]
  = iZ\}$.  Applying Lemma \ref{sub:sl2-lemma} we see that
  $\ad(X_\FF):\plie_{\overline{\Q}} \to \vlie_{\overline{\Q}}$ is
  surjective, from which it follows that $\ad(X_\FF):\plie_\FF \to
  \vlie_\FF$ is surjective.
  
  Proposition \ref{sub:spaltenstein}(1) shows that $\dim_k
  \clie_\plie(X) = \dim_k \clie_\glie(X)$. If we regard $\ad(X_\FF)$
  and $\ad(X)$ as endomorphisms respectively of $\plie_\FF$ and of
  $\plie$, then (2) of that proposition yields $\dim_\FF \ker
  \ad(X_\FF) = \dim_k \ker \ad(X)$; thus we apply Lemma
  \ref{sub:order-generalities} to conclude that
  $\ad(X)^{n(P)-1}(\plie) \not = 0$ as desired.
\end{proof}

\section{The exponent of $\V$}
\label{sec:Group}

Recall that we have fixed $G$ a reductive group in good
characteristic, and $P$ a distinguished parabolic subgroup with Levi
decomposition $P = \Levi \V$.

\subsection{A Theorem of Springer}
\label{sub:springer-theorem}

We recall the following important result.
\begin{prop}
  Let $G$ be quasisimple, and assume that $p$ is a good prime for $G$.
  If the root system is of type $A$, assume that the isogeny $G_{sc}
  \to G$ is separable.  Then there is a $B$-equivariant isomorphism of
  varieties $\e:\ulie \to U$ which extends to a $G$-equivariant
  isomorphism of varieties $\e$ between the nilpotent variety of
  $\glie$ and the unipotent variety of $G$.
\end{prop}
A version of the proposition was first proved by Springer
\cite{springer-unipotent-iso}, though with the slightly weaker
conclusion that $\e$ is a homeomorphism.  We refer the reader to the
discussion of this result in \cite[6.20/1]{hum-conjugacy}.

In view of the equivariance
property, it is clear that $\e$ restricts to a $P$-isomorphism
$\vlie \to \V$, where $\vlie \subset \ulie$ is the Lie algebra of $\V$.
If $X \in \vlie$ is a Richardson element, we let $R(X)$ denote
the unipotent radical of the centralizer of $X$ in $G$; note that
$R(X)$ coincides with the centralizer in $\V$ of $X$.

\begin{cor}
  Assume that $G$ satisfies the hypothesis of the proposition, and let
  $X \in \vlie$ be a Richardson nilpotent element. 
\begin{enumerate}
\item Let $Z$ be the center of $R(X)$. Then $Z$ is a connected group,
  and $X$ is tangent to $Z$; i.e. $X \in \Lie(Z)$.
\item The exponent of the connected, Abelian, unipotent group $Z$ is
  $\ge p^n$ where $n$ is the $p$-nilpotence degree of $X$.
\end{enumerate}
\end{cor}

\begin{proof}
  Put
  \begin{equation*}
    \wlie = \{Y \in \vlie \mid \Ad(u)Y = Y \ 
    \text{for all $u \in R(X)$}\},
  \end{equation*}
  and 
  \begin{equation*}
    W = \{y \in \V \mid u^{-1}yu = y \ 
    \text{for all $u \in R(X)$}\}.
  \end{equation*}
  Then $\wlie$ is a linear subspace of $\vlie$, and $\e(\wlie) = W$.
  This shows that $W$ is a connected subgroup of $\V$.
  
  Denoting by $\zlie$ the Lie algebra of $Z$, we have $\zlie \subset
  \wlie$ since any $u$ in $R(X)$ centralizes $X$ by definition.
  Similarly, we have $Z \subset W$. On the other hand, if $w \in W$,
  then $w^{-1}\e(X)w = \e(X)$ since $\e(X) \in R(X)$; this shows that
  $\Ad(w)X = X$ hence $w \in R(X)$. Since $w$ commutes with each
  element of $R(X)$, we deduce that $w \in Z$ hence $Z=W$. This shows
  that $Z$ is connected, and that $\dim \zlie = \dim Z = \dim W = \dim
  \wlie$ so that $\zlie = \wlie$. Since $X \in \wlie$, we get $X \in
  \zlie$ and (1) follows.

  To see the second assertion of the corollary, one applies
  Proposition \ref{sub:connect-abelian}.
\end{proof}

\subsection{} 
We can now prove the \emph{Order Formula} for unipotent elements
originally obtained by D. Testerman in \cite{testerman}.

\label{sub:order}
\begin{theorem}
  \emph{(Order Formula)} Let $m \ge 1$ be the minimal integer with
  $p^m \ge n(P)$, as in Theorem \ref{sub:p-exponent}. Then the
  exponent of the group $\V$ is $p^m$, and the order of a
  Richardson element of $\V$ is $p^m$.
\end{theorem}

\begin{proof}
  In view of Lemma \ref{sub:parabolic}, we may suppose that $G$
  satisfies the hypothesis of Proposition \ref{sub:spaltenstein}
  and of Proposition \ref{sub:springer-theorem}.
  
  Let $Z \le \V$ be as in Corollary \ref{sub:springer-theorem}.  Then
  that corollary together with Theorem \ref{sub:p-exponent} imply that
  \begin{equation*}
    p^m \le \text{exponent of $Z$} \le \text{exponent of $\V$} \le p^m
  \end{equation*}
  whence equality holds.
\end{proof}

\begin{example}
  Let $G$ be a group of type $G_2$, and let $p \ge 5$, so that $p$ is
  good for $G$. There are two distinguished orbits of nilpotent (and
  unipotent) elements; these are usually labelled $G_2$ (for the
  regular class) and $G_2(a_1)$ (for the subregular class). We choose
  a fixed Borel subgroup $B$; a Richardson element for $B$
  represents the regular class. If $P$ is a minimal parabolic subgroup
  containing $B$, a Richardson element for $P$ represents the
  subregular class.  Moreover, $P$ is distinguished only if
  $P=P_\alpha$ where $\alpha$ is the short simple root.  We have $n(B)
  = h = 6$ and $n(P_\alpha) = 4$. A subregular unipotent element
  always has order $p$ and a subregular nilpotent element is
  $p$-nilpotent. A regular unipotent element has order $p$ unless
  $p=5$ in which case it has order 25; likewise, a regular nilpotent
  element is $p$-nilpotent unless $p=5$ in which case it has
  $5$-nilpotence degree 2.
\end{example}

\begin{rem}
  One can compute the order of an arbitrary unipotent $u$ by finding a
  pair $\Levi,Q$ where $\Levi \le G$ is a Levi subgroup containing $u$
  and $Q$ is a distinguished parabolic subgroup of $\Levi$ for which
  $u$ is a Richardson element.  Similar remarks hold for an arbitrary
  nilpotent $X$.
\end{rem}

\section{Exponentials in Linear Algebraic Groups}
\label{sec:exponentials}

\subsection{Exponentials in characteristic 0.}
\label{sub:exp-char-0}

Let $\FF$ be a field of characteristic 0, and let $G$ be a linear
algebraic group defined over $\FF$.  For any nilpotent element $X \in
\glie_\FF$ (the $\FF$-form of $\glie = \Lie(G)$) and any rational
representation $(\rho,V)$ of $G$, $d\rho(X)$ is nilpotent so one may
define a homomorphism of algebraic groups
\begin{equation*}
  \e_{X,V}:\G_a \to \GL(V) \quad \text{via}\ t \mapsto \exp(d\rho(tX))
\end{equation*}
by the usual formula. If $(\rho,V)$ is defined over $\FF$, then so is
$\e_{X,V}$.

\begin{prop}
  There is a unique homomorphism of algebraic groups $\e_X:\G_a \to G$
  such that $\e_{X,V} = \rho \circ \e_X$ for all rational
  representations $(\rho,V)$ of $G$. The homomorphism $\e_X$ is
  defined over $\FF$.
\end{prop}

\begin{proof}
  Let $(\rho,V)$ be a faithful $\FF$-representation of $G$. Since the
  image of $d\e_{X,V}$ is a subalgebra of $d\rho(\glie_\FF) \subset
  \gl(V)$, one gets $\e_{X,V}(\G_a) \le \rho(G) \le \GL(V)$ by
  \cite[Cor.  6.12]{Bor1}; thus we get a morphism $\e_X:\G_a \to G$
  defined over $\FF$ and satisfying $\rho \circ \e_X = \e_{X,V}$.  A
  second application of the result of \emph{loc. cit.} shows that
  $\rho'\circ\e_X=\e_{X,V'}$ for any rational representation
  $(\rho',V')$. Unicity of $\e_X$ is clear.
\end{proof}

\subsection{Exponentials over integers}
\label{sec:integral-exp}

Let $A$ be a Dedekind domain, with field of fractions $\FF$.  We
suppose that $\FF$ has characteristic 0. [Important examples of $A$
for us are: the rational integers $\Z$, the ring of integers in a
number field $\FF$, a localization of a Dedekind domain at a prime
ideal].

Let $G_\FF$ denote a linear algebraic group over $\FF$, with a
faithful $\FF$-representation $(\rho,V)$. Fix an $A$-lattice $\LL
\subset V$; thus $\LL$ is a finitely generated $A$-module containing
an $\FF$-basis of $V$. If we localize at a maximal ideal $\mm$,
then $\LL_\mm$ is a free $A_\mm$ module; thus $\LL$ is
$A$-projective.

Let $J \normal \FF[\GL(V)]$ denote the ideal defining $G_\FF$
(more precisely: defining $\rho(G_\FF)$). The choice of $A$-lattice
$\LL$ determines the integral form $A[\GL(\LL)] \subset \FF[\GL(V)]$;
let $J_A = J \cap A[\GL(\LL)]$.

We make the following assumption:
\begin{equation}
\label{eq:gp-scheme-hyp}
  \text{The $A$-algebra $B = A[\GL(\LL)]/J_A$ represents a group scheme $G_A$.}
\end{equation}

\begin{prop}
  \label{prop:gp-scheme-factorize}
  Let $H_A$ be an affine group scheme over $A$ such that $A[H_A]$ is
  free as an $A$-module. Let $\phi:H_A \to \GL(\LL)$ be a
  homomorphism of group schemes. If the base-changed map $\phi_\FF$
  determines a morphism $H_\FF \to G_\FF$, then $\phi$ determines
  a morphism $H_A \to G_A$.
\end{prop}

\begin{proof}
  $\phi$ is determined by its comorphism $\phi^*:A[\GL(\LL)] \to
  A[H_A]$; the proposition will follow if we show that $\phi^*(J_A) =
  0$. Note that $\FF[\GL(V_\FF)] \iso A[\GL(\LL)] \tensor_A \FF$ and
  $\FF[H_\FF] = A[H_A] \tensor_A \FF$ (the latter by definition). The
  comorphism $\phi^*_\FF$ is then $\phi^* \tensor 1$.  The hypothesis
  implies that $\phi_\FF^*(J) = 0$; since $A[H_A]$ is free as an
  $A$-module, the natural map $A[H_A] \to \FF[H_\FF]$ is injective,
  and it then follows that $\phi^*(J_A) = 0$ as desired.
\end{proof}

\begin{cor}
  Let $X \in \glie_\FF$ be nilpotent, and suppose that $d\rho(X) \in
  \End_A(\LL)$.
  \begin{enumerate}
  \item[(a)] Suppose that
    $\exp(d\rho(X))\LL \subseteq \LL$. Then the exponential
    homomorphism 
    \begin{equation*}
      t \mapsto \exp(d\rho(tX)):\G_{a,A} \to \GL(\LL)
    \end{equation*}
    determines a homomorphism of group schemes $\G_{a,A} \to G_A$ over
    $A$.
  \item[(b)] Let $p$ be a rational prime, and let $\mm \normal A$ be a
    maximal ideal for which $A/\mm$ has characteristic $p$.
    Suppose that $d\rho(X)^{p^i} \in d\rho(\glie_\FF)$ for $0 \le i <
    n$, and that $d\rho(X)^{p^n} = 0$. Then the Artin-Hasse
    exponential map (see \ref{sub:Artin-Hasse-series})
    \begin{equation*}
      \vec{t} \mapsto
      E_{d\rho(X)}(\vec{t}):\Witt_{n,A_\mm} \to \GL(\LL_\mm)
    \end{equation*}
    determines a homomorphism of group schemes $\Witt_{n,A_\mm} \to
    G_{A_\mm}$ over $A_\mm$.
  \end{enumerate}
\end{cor}

\begin{proof}
  Note first that the coordinate algebras $A[\G_a]$ and $A[\Witt_n]$
  are free as $A$ modules.  (a) follows immediately from the proposition
  combined with Proposition \ref{sub:exp-char-0}. For (b), Lemma
  \ref{sub:Artin-Hasse-linear}(1) combined with Proposition
  \ref{sub:exp-char-0} shows that $E_{d\rho(X)}$ determines on base
  change a morphism $\Witt_{n,\FF} \to G_\FF$; the result then follows
  from the proposition.
\end{proof}

\subsection{Nilpotent orbits and fields of definition}
\label{sub:fields-o-def}

If $k \subset k'$ are two algebraically closed fields, then an
algebraic group $G_k$ over $k$ determines by extension of scalars an
algebraic group $G_{k'}$ over $k'$.  Suppose that $G_k$ acts on an
affine variety $V_k$; $G_{k'}$ also acts on $V_{k'}$.

The following result is attributed to P. Deligne in the introduction
to G. Lusztig's paper ``On the finiteness of the number of unipotent
classes," [Invent. Math. 34 (1976)].  A proof due to R.
Guralnick can be found in \cite[Prop
1.1]{GLMS-finite-orbits}.
\begin{prop}
  Suppose that $G_k$ has finitely many orbits on $V_k$. Then $G_{k'}$
  has finitely many orbits on $V_{k'}$, and each $G_{k'}$ orbit has a
  $k$-rational point. In particular, the number of $G_k$ orbits on
  $V_k$ is the same as the number of $G_{k'}$ orbits on $V_{k'}$.
\end{prop}

Richardson's theorem \cite[Theorem 3.10]{hum-conjugacy} (together with
case-by-case analysis for bad primes -- see the discussion in
\emph{loc. cit.} Theorem 6.19) shows that a reductive group has
finitely many orbits on its nilpotent variety, so we obtain:
\begin{cor}
  If $G$ is a reductive group over an algebraically closed field $k$,
  then each nilpotent orbit of $G$ in $\glie$ contains a point which
  is rational over the algebraic closure of the prime field in $k$.
\end{cor}

\begin{rem}
  When $p=0$ or is sufficiently large for the reductive group $G$, it
  follows from \cite[Theorem III.4.29]{springer-steinberg} that each
  nilpotent orbit has a point rational over the prime field.  We don't
  need this fact.
\end{rem}

\subsection{Exponentials and Chevalley groups}
\label{sub:admiss-group}

We have already mentioned in \ref{sub:spaltenstein} the existence of a
split reductive group scheme $G_\Z$ over $\Z$ from which $G$ arises by
base change; we now need more precise information about $G_\Z$.

We suppose $G$ to be a semisimple group over $k$. Then $G$ is
(isomorphic with) a Chevalley group; we recall some of the ideas
behind this construction (for which the reader may find full details
in \cite{Steinberg}).

Let $\glie_\Q$ denote a split simple Lie algebra over $\Q$ with the
same root system as $G$. For a suitable finite dimensional
$\Q$-representation $(d\rho,V)$ of $\glie_\Q$, and a $\Z$-lattice $\LL
\subset V$ invariant by Kostant's $\Z$-form of the enveloping algebra
of $\glie_\Q$, one ``exponentiates'' the action of Chevalley basis
elements on $\LL$ to obtain Chevalley groups with the following
properties:
\begin{itemize}
\item Over $\Q$, one gets a closed subgroup $G_\Q \le \GL(V)$ defined
  and split semisimple over $\Q$. The root datum of $G_\Q$ is the same
  as that of $G$. Moreover, the $\Q$-Lie algebra of $G_\Q$ is
  $\glie_\Q$.
\item In characteristic $p>0$, one gets a closed subgroup $G_{\F_p}
  \le \GL(\LL/p\LL)$, defined and split semisimple over $\F_p$, which
  is isomorphic over $k$ with the original group $G$.
\end{itemize}

Since there should be no danger of confusion, we will write $(\rho,V)$
for the representation of $G_\Q$ on $V$, and $(\rho,\LL/p\LL)$ for the
representation of $G_{\F_p}$ on $\LL/p\LL$.

As in \ref{sec:integral-exp}, let $J \normal \Q[\GL(V)]$ be the ideal
defining the $\Q$-variety $G$, and put $J_\Z = J \cap \Z[\GL(\LL)]$
(where $\Z[\GL(\LL)]$ is regarded as a subring of $\Q[\GL(V)]$ in the
obvious way).  Let $B$ be the $\Z$-algebra $\Z[\GL(\LL)]/J_\Z$.
\begin{lem}
  The $\Z$-algebra $B$ represents a split semisimple group scheme
  $G_\Z$ over $\Z$. One has $B \tensor_\Z \F_p = \F_p[G_{\F_p}]$ and
  $B \tensor_\Z \Q = \Q[G_\Q]$, or equivalently, $G_\Q$ and $G_{\F_p}$
  are obtained by base change from $G_\Z$.
\end{lem}
\begin{proof}
  This is proved in \cite[\S3.4 and
  \S4]{borel70:_proper_linear_repres_cheval_group}.
\end{proof}

Let now $A$ be a Dedekind domain with field of fractions $\FF$ as in
\ref{sec:integral-exp}.  We suppose that $\FF$ is a finite extension
of $\Q$ (so $\FF$ is a number field). We regard $\FF$ as a subfield of a
fixed algebraic closure $\overline{\Q}$ of $\Q$. Let $G_A$ be
the group scheme obtained from $G_\Z$ by base change; thus $G_A$ is
represented by
\begin{equation}
  B \tensor_\Z A = A[\GL(\LL \tensor_\Z A)]/J_A
\end{equation}
[where $J_A$ is the ideal $\FF J \cap A[\GL(\LL_A)]$ with the
intersection occurring in $\FF[\GL(V_\FF)]$.]

The choice of a Chevalley basis of $\glie_\Q$ entails the choice of a
triangular decomposition $\glie_\Q = \ulie^-_\Q \oplus \hlie_\Q \oplus
\ulie_\Q$ where $\hlie_\Q$ is a maximal toral subalgebra; the
construction of $G_\Q$ yields also a maximal torus $T_\Q \le G_\Q$
with $\Lie(T_\Q) = \hlie_\Q$. Moreover, we have the decomposition
$\glie_\Z = \ulie^-_\Z \oplus \hlie_\Z \oplus \ulie_\Z$, where e.g.
$\ulie_\Z$ is the $\Z$-span of the Chevalley basis elements
which correspond to the positive roots.

Let $\phi = \sum_{\alpha > 0} \alpha^\vee$ regarded as a cocharacter
for the torus $T_\Q$.  To the $\glie_\Q$ module $V$, we associate the
integer
\begin{equation*}
  n(V) = \max\{\langle \lambda,\phi \rangle \mid \lambda \in X^*(T_\Q),
  \ V_\lambda \not = 0\}.
\end{equation*}
\begin{prop}
  \begin{enumerate}
  \item If $X \in \glie_{\overline{\Q}}$ is nilpotent, then
    $d\rho(X)^{n(V)+1} = 0$.
  \item Suppose $X \in \glie_\FF$ is nilpotent and satisfies
    $d\rho(X)\LL_A \subseteq \LL_A$. If $n(V)!$ is invertible in $A$
    (e.g. if $A$ is local with residue field of characteristic $p >
    n(P)$), then $\exp(d\rho(X))$ leaves $\LL_A$
    invariant and $t \mapsto \exp(d\rho(tX))$ defines a morphism of
    group schemes $\G_{a,A} \to G_A$.
  \end{enumerate}
\end{prop}

\begin{proof}
  (1). We may suppose that $X \in \ulie_{\overline{\Q}}$. The
  co-character $\phi$ induces a grading on $V_{\overline{\Q}}$ by $V_i
  = \{v \in V_{\overline{\Q}} \mid \rho(\phi(t))v = t^iv \ \forall t
  \in \overline{\Q}\}$; evidently $d\rho(X)$ acts as a sum of
  homogeneous terms of positive and even degree. Since the Weyl group
  of $\glie_\Q$ permutes the weights of $V$, it follows that
  $V_{\overline{\Q}} = \bigoplus_{-n(V) \le i \le n(V)} V_i$, and (1)
  is then immediate.  [One can alternately argue that the graded
  components of $V_{\overline{\Q}}$ are the weight spaces for the
  action of an $\lie{sl}_2$-subalgebra containing a regular nilpotent
  element, so that $n(V)$ is the highest weight. Since the regular
  nilpotent elements are dense in the nilpotent variety, (1) follows;
  moreover, this argument shows that $d\rho(X)^{n(V)} \not = 0$ when
  $X$ is regular.]
  
  (2) Since $i!$ is invertible in $A$ for each $0 \le i \le n(V)$ and
  since $d\rho(X)^{n(V)+1} = 0$ by (1), it follows that
  $\exp(d\rho(X))$ leaves $\LL_A$ invariant.  The result now holds by
  Corollary \ref{sec:integral-exp}(a).
\end{proof}

\begin{cor}  
  Suppose that $n(V) < p$ and that $X \in \glie_k$ is nilpotent.  Then
  $X^{[p]} = 0$, and the (truncated) exponential $t \mapsto
  \exp(d\rho(tX))$ defines a morphism of algebraic groups $\G_{a,k}
  \to G = G_k$.
\end{cor}

\begin{proof}
  In view of the results of \ref{sub:fields-o-def}, we may suppose
  that $k$ is an algebraic closure of the finite field $\F_p$.
  
  The image of $\ulie_\Z$ in $\glie_{\F_p}$ is the $\F_p$-Lie algebra
  of the unipotent radical of a Borel subgroup.  Since $k$ is an
  algebraic closure of $\F_p$, there is a finite extension $\FF$ of
  $\Q$ and a valuation ring $A$ in $\FF$ whose residue field $A/\mm =
  l \subset k$ has the property that $\Ad(g)X \in \ulie_l$ for a
  suitable $g \in G(l)$. Thus, we may suppose $X \in \ulie_l$.  Since
  $\ulie_l = \ulie_A/\mm \cdot \ulie_A$, we may choose a lift $\tilde
  X \in \ulie_A$ of $X$. Since each element of $\ulie_\FF$ is
  nilpotent, part (1) of the previous proposition now shows that
  $d\rho(\tilde X)^p = 0$, hence also $d\rho(X)^p = 0$ which implies
  $X^{[p]}= 0$.
  
  Since $n(V)$ is invertible in $A$, part (2) of the previous
  proposition shows that the exponential map determines a morphism of
  group schemes $\e:\G_{a,A} \to G_A$ over $A$.  The condition
  $d\rho(\tilde X)^p = 0$ implies that $\rho \circ \e$ has degree $<
  p$ when regarded as a morphism of $\FF$-varieties $\G_{a,\FF} =
  \Aff^1 \to \GL(V)$.  Denoting by $\overline{\e}:\G_{a,k} \to G_k$
  the morphism obtained by base change, it follows that also $\rho
  \circ \overline{\e}$ has degree $< p$.
  
  The differential $d\e:\FF = \Lie(\G_{a,\FF}) \to \glie_\FF$
  satisfies $d\e(1) = \tilde X$. It follows that $d(\rho \circ
  \overline{e})(1) = d\rho(X)$.
  
  Now, the exponential through $d\rho(X)$ is the unique homomorphism
  $\G_{a,k} \to \GL(\LL_A \tensor_A k)$ with degree $<p$ and whose
  differential at 1 is $d\rho(X)$. Thus, $\overline{\e}$ coincides
    with the truncated exponential, and the corollary is proved.
\end{proof}

\begin{rems}
  \begin{enumerate}
  \item Suppose that $G$ is of adjoint type (i.e. that the character
    group of a maximal torus is spanned over $\Z$ by the roots). Then
    $G$ arises as a Chevalley group where we may take the adjoint
    module $(\ad,\glie_\Q)$ for the $\glie_\Q$ module $(d\rho,V)$.  In
    this case, one has $n(\glie_\Q) = 2h-2$ where $h$ is the Coxeter
    number of the root system of $G$ (see \cite[Prop.
    5.5.2]{Carter1}). 
    
    The reader should compare the result of the Corollary in this case
    with \cite[Prop. 5.5.5(iv)]{Carter1}, where a similar conclusion
    is asserted, but with the weaker bound $p > 3h-3$.  We emphasize
    that the argument doesn't depend on the Bala-Carter theorem.  We
    have used the finiteness of nilpotent orbits in order to know that
    every nilpotent orbit has a rational point over the algebraic
    closure of a finite field; as noted before, in good
    characteristic that finiteness is a result of Richardson and is
    independent of the classification of nilpotent orbits.

  \item We will be most interested in applying the corollary in case
    $G$ is quasisimple with exceptional root system.  We list here
    some data for these root systems, including the minimal
    dimensional non-trivial $\glie_\Q$ module $V_{\text{min}}$.  The
    module $V_{\text{min}}$ is a simple module $L_\Q(\lambda)$ with
    highest $\lambda$; we describe $\lambda$ in terms of the
    fundamental dominant weights with the labelling as in the tables
    in \cite[Planche V-IX]{BouLie456}.  Those tables may be used to
    compute the indicated values of $n(V_{\text{min}}) = \langle
    \lambda,\varphi \rangle$.
    \begin{equation*}
      \begin{array}[t]{lclll}
        R & 2h-2 & V_{\text{min}} & n(V_{\text{min}})  \\
        \hline \hline
        G_2  & 10 & L_\Q(\varpi_1) & 6 \\
        F_4  & 22 & L_\Q(\varpi_4) & 16 \\
        E_6  & 22 & L_\Q(\varpi_1) & 16 \\
        E_7  & 34 & L_\Q(\varpi_7) & 27 \\
        E_8  & 58 & L_\Q(\varpi_8) & 58 \\
      \end{array}
    \end{equation*}
    
  \item The technique used in proof of the Theorem is similar to that
    used in \cite{testerman}. Let $Y \in \glie_\Q$ with $d\rho(Y) \in
    \End_{\Zp}(\LL)$. We remark that \cite[Lemma
    1.4]{testerman} gives moreover a condition under which
    $\exp(d\rho(Y))$ leaves invariant the lattice $\LL$ even when
    $d\rho(Y)^p \not = 0.$ We mention also a different condition:
    namely, $(*)$ if $d\rho(Y)^p\LL \subseteq p\LL$ and
    $d\rho(Y)^{p^2} = 0$ then $\exp(d\rho(Y))$ leaves $\LL_{\Z_{(p)}}$
    invariant.  This follows from the formula for $\nu_p(i!)$ which
    may be found in \cite[Ch. I, Exerc. 13c]{Koblitz-p-adic}.
    
    The condition $d\rho(Y)^p\LL\subseteq p\LL$ means that the image
    $\overline{Y}$ of $Y$ in $\glie_k$ is $p$-nilpotent; one may
    argue, as in the corollary, that under the condition $(*)$ there
    is a homomorphism $\G_a \to G$ over $k$ obtained by base change
    from the exponential over $\Zp$. However, the homomorphism over
    $\Zp$ need not have degree $<p$; thus the map obtained by
    reduction modulo $p$ need not coincide with the truncated
    exponential of $\overline{Y}$ in $\GL(\LL/p\LL)$.
  \end{enumerate}
\end{rems}

\subsection{Classical groups and the Artin-Hasse exponential.}
\label{sub:classical-witt}

Let $V$ be a $\Q$-vector space with a bilinear form $\varphi$. 
Let $G_\Q$ be the stabilizer in $\SL(V)$ of $\varphi$. We
assume that one of the following three statements holds: 
\begin{enumerate}
\item[CG1.] $\varphi = 0$, so that $G_\Q = \SL(V)$
\item[CG2.] $\varphi$ is non-degenerate and alternating, so that $G_\Q
  = \text{Sp}(V,\varphi)$,
\item[CG3.] $\varphi$ is non-degenerate and symmetric, so that $G_\Q =
  \text{SO}(V,\varphi)$. Moreover, if $\dim_\Q V$ is written as $2r +
  \epsilon$ with $\epsilon \in \{0,1\}$, then $V$ contains a totally
  singular $\Q$-subspace of dimension $r$ (so $\varphi$ has maximal
  Witt index, or is a \emph{split form}).
\end{enumerate}
In each case, $G_\Q$ is a connected, quasisimple $\Q$-split group.

The Lie algebra $\glie_\Q$ of $G_\Q$ is split simple, and we may carry
out the ``Chevalley group'' constructions of \ref{sub:admiss-group}
for $\glie_\Q$ with respect to its natural representation $(\nu,V)$.
Fix a lattice $\LL \subset V$ invariant by $\mathcal{U}_\Z$. Over
$\Q$, the group constructed in this way identifies with the original
group $G_\Q$; this follows from \cite{Ree-chev-gp}.  For each prime
$p$, we get also a quasisimple $\F_p$-split algebraic group $G_{\F_p}$
with the same root datum as $G_\Q$.

The formal character of the $G_{\F_p}$-representation $(\nu,\LL/p\LL)$
coincides with the formal character of the $G_\Q$-representation
$(\nu,V)$; moreover, it is well known that $G_{\F_p}$ has an
irreducible representation with that formal character (provided $p\not
= 2$ in case CG3 with $\epsilon =1$).  Thus $\LL/p\LL$ is irreducible
for $G_{\F_p}$ (with this restriction on $p$).

Replacing $\varphi$ by a suitable integral multiple, we may suppose in
case CG2 or CG3 that $\varphi(\LL,\LL) = \Z$; note that the group
$G_\Q$ and Lie algebra $\glie_\Q$ are unchanged by this replacement.
For each prime $p$ (with the restriction on $p$ of the previous
paragraph), $\varphi$ induces a non-0 bilinear form
$\overline{\varphi}$ on $\LL/p\LL$; since that module is irreducible,
$\overline{\varphi}$ must be non-degenerate.  It then follows from
\cite{Ree-chev-gp} that $G_{\F_p}$ coincides with the stabilizer in
$\SL(\LL/p\LL)$ of $\overline{\varphi}$ so long as $p\not = 2$ in case
CG3 (for either value of $\epsilon$).

Let $B = \Z[\GL(\LL)]/J_\Z$ as in \ref{sub:admiss-group}.  Then $B$
represents a group scheme over $\Z$, and for each prime integer $p$,
it follows from Lemma \ref{sub:admiss-group} that $B \tensor_\Z \F_p$
is the coordinate ring of $G_{\F_p}$.  

\begin{lem}
  Let $n$ be a positive integer, and suppose that $n$ is odd in case
  (ii) or (iii). Let $\FF$ be an extension field of $\Q$. If $X \in
  \glie_{\FF}$, then $d\nu(X)^n \in d\nu(\glie_\FF)$.
\end{lem}

\begin{proof}
  We have for each $v,w \in V_\FF$:
  \begin{equation*}
    \varphi(d\nu(X)^nv,w) = (-1)^n\varphi(v,d\nu(X)^nw),
  \end{equation*}
  whence the result.
\end{proof}

We now deduce in this case a
recent result of R. Proud:
\begin{prop}
  Suppose that $p$ is good for $G_k$ (i.e. that $p \not = 2$ in case
  CG2 or CG3). For each nilpotent $X \in \glie$ with $p$-nilpotence
  degree $n$, the Artin-Hasse exponential defines an injective
  morphism of algebraic groups $E_X:\Witt_n \to G_k$.  Thus each
  unipotent element of $G$ lies in a closed subgroup isomorphic with
  some $\Witt_n$. If $l$ is a subfield of $k$ and $X \in \glie_l$ then
  $E_X$ is defined over $l$.
\end{prop}

\begin{proof}
  Using the results of \ref{sub:fields-o-def}, we may suppose
  that $k$ is an algebraic closure of the finite field $\F_p$.  As in
  the proof of Corollary \ref{sub:admiss-group}, we may find a number
  field $\FF$ with valuation ring $A$ and residue field $l = A/\mm$
  for which $X$ lies in $\ulie_l$ (where $\ulie_\Z$ is the $\Z$-span
  of suitable Chevalley basis elements, as before). Thus we may choose
  a lift $\tilde X \in \ulie_A$ of $X \in \ulie_l = \ulie_A /\mm
  \ulie_A$.
  
  Corollary \ref{sec:integral-exp}(b) now yields a homomorphism of
  group schemes $E_{\tilde X}:\Witt_{m,A} \to G_A$ given by the
  Artin-Hasse exponential, where $m \ge n$ is the nilpotence degree of
  $d\nu(\tilde X)$.
  
  We get then by base change a homomorphism $\Witt_{m,k} \to G_k$ over
  $k$; note that by the formula defining $E_{\tilde X}$ in
  \ref{sub:Artin-Hasse-linear}, this base-changed homomorphism
  vanishes on the subgroup $K = \{(0,\dots,0,t_n,\dots,t_{m-1})\} \le
  \Witt_{m,k}$, and coincides with the homomorphism
  $E_{d\nu(X)}:\Witt_{m,k}/K = \Witt_{n,k} \to \GL(V)$. It follows
  that $E_{d\nu(X)}$ takes values in $G_k$ (hence has rights to be
  called $E_X$), and is injective, as claimed.
  
  It is clear that the partition of $\dim V$ determined by the Jordan
  block sizes of the unipotent element $\nu(E_X(1))$ on the natural
  module $V$ is the same as the partition of $d\nu(X)$. In the cases
  CG1, CG2 and CG3 with $\epsilon = 1$, the unipotent classes of $G$
  and nilpotent classes of $\glie$ are classified by these partitions
  (see \cite[7.11]{hum-conjugacy}), so we get the claim on unipotent
  elements in these cases. In case CG3 with $\dim V$ even, let
  $G'=O(\LL \tensor_A k)$ denote the full orthogonal group; thus $G$
  is the identity component of $G'$, and has index 2. The unipotent
  elements of $G'$ all lie in $G$, and $\Lie(G) = \Lie(G')$. Again by
  \cite[7.11]{hum-conjugacy}, the unipotent and nilpotent classes of
  $G'$ are classified by partition, so it is clear that each unipotent
  element of $G'$ lies in a suitable Witt-vector subgroup. But any
  such subgroup, being connected, must lie in $G$.
  
  The rationality assertion is clear.
\end{proof}

R. Proud has proved the proposition for \emph{all} quasisimple groups
$G$, not only classical groups; see \cite{Proud-Witt}.  His techniques
for classical $G$ are different from those used here.

\section{The exponential isomorphism for $p \ge n(P)$}
\label{sec:exponential-big-p}

We describe in this section an argument due to Serre \cite[\S
2.2]{serre:PSLp} that will be used below.  This argument is also
used in some recent work of Gary Seitz \cite[\S 5]{seitz-unipotent}.

Suppose $B_\Z$ is a Borel subgroup of the split reductive group $G_\Z$
over $\Z$, and let $B \le G$ be the corresponding groups over $k$.
For $\alpha \in R$, let $\phi_\alpha$ be an isomorphism over $\Z$
between $\G_a$ and the root subgroup $U_\alpha \le B$.

Any standard parabolic subgroup $B \le P \le G$ is defined over $\Z$.
If $\V$ denotes the unipotent radical of $P$, then the $\phi_\alpha$
define an isomorphism $\prod_{\alpha \in R \setminus R_I}
U_{\alpha,\Z} \to \V_\Z$ of schemes over $\Z$ (where $I \subseteq S$
defines the parabolic subgroup $P$ as in \ref{sub:parabolic}).  For
any $\Z$-algebra $\Lambda$, a point $u$ of $\V$ over $\Z$ may thus be
written uniquely as $u = \prod_{\alpha \in R \setminus R_I}
\phi_\alpha(t_\alpha)$ with $t_\alpha \in \Lambda$; thus the
$t_\alpha$ form a system of coordinates for $\V$ over $\Z$.

The Lie algebra $\vlie_\Z$ is the $\Z$-span of the $e_\alpha =
d\phi_\alpha(1)$ for $\alpha \in R^+ \setminus R^+_I$.
The nilpotence degree $n = n(P)$ of $\V_\Q$ (and of $\vlie_\Q$)
is given by the formula in \ref{sub:nilp-class}; we will
work with the ring $A = \Z[1/(n-1)!]$.

Proposition \ref{sub:exp-char-0} implies that the exponential map
defines a morphism of varieties $\e:\vlie_\Q \to \V_\Q$ over $\Q$.
Similarly, the logarithm yields a morphism of varieties $\V_\Q \to
\vlie_\Q$, so that $\e$ is an isomorphism.  For each $\Q$-algebra
$\Lambda$, each $\Lambda$-point $u$ of $\V$ may be written uniquely as
$u = \e(\sum_{\alpha \in R \setminus R_I} u_\alpha e_\alpha)$ for
$u_\alpha \in \Lambda$. Thus the $u_\alpha$ form a system of
coordinates for $\V$ over $\Q$.

Since $\vlie_\Q$ is a nilpotent Lie algebra, it may be regarded as an
algebraic group over $\Q$ via the Hausdorff series (compare
\cite[2.2]{serre:PSLp} for the case $P = B$, and see \cite[Ch. 2,
\S6]{BouLie123}).  Moreover, it follows from \cite[ch. 2 \S 6.4 Theorem
2]{BouLie123} that the operation in $\vlie_\Q$ is defined over $A$,
so that $\vlie_A$ is an affine group scheme over $A$.

\begin{theorem}
  The exponential map $\e$ defines a $P_A$-equivariant isomorphism of
  group schemes $\vlie_A \to \V_A$.
\end{theorem}

\begin{proof}
  The essential point is proved in \cite[\S2.2, Prop. 1]{serre:PSLp}
  for $P=B$; the generalization to $P$ is immediate. One observes as
  in \emph{loc. cit.} that
  \begin{equation*}
    u_\alpha = t_\alpha + P_\alpha((t_\beta)_{\beta < \alpha})
  \end{equation*}
  where $P_\alpha$ is a polynomial with coefficients in $A =
  \Z[1/(n-1)!]$ in the $t_\beta$ with $\beta < \alpha$. It follows
  that $\e$ is an isomorphism over $A$ (see \cite[\S2.2, Rem.
  2]{serre:PSLp}). The equivariance assertion is clear.
\end{proof}

\begin{example}
  Let $G = \SP_4(\Q)$, so that $R$ is of type $C_2$, and let $P = B$.
  Recall that $R^+ = \{\alpha,\beta,\alpha+\beta,2\alpha+ \beta\}$.
  It is straightforward to check that
  \begin{equation*}
      \e(X ) =
    \phi_\alpha(a)\phi_\beta(b) \phi_{\alpha +
      \beta}\left(c + \frac{ab}{2}\right) 
    \phi_{2\alpha +
      \beta}\left(d - bc
      - \frac{2ab^2}{3}\right)
  \end{equation*}
  for $X = a e_\alpha + b e_\beta + c e_{\alpha + \beta} + d
  e_{2\alpha + \beta}$.  It is then clear that $\exp$ is
  an isomorphism over $A = \Z[1/6]$ (note that $h-1=n(B)-1 = 3$).
\end{example}

If $p \ge n(P)$, then the field $k$ is an $A$
algebra (in a unique way), and the above result yields the following:
\begin{cor}
  \label{cor:exp-iso} Suppose that $p \ge n(P)$.
  \begin{enumerate}
  \item $\e$ is a $P$-equivariant isomorphism of $k$-varieties
  $\vlie \to \V$.  
\item If $X,Y \in \vlie$ satisfy $[X,Y] = 0$, then $\e(X)$ and
  $\e(Y)$ commute.
  \end{enumerate}
\end{cor}

\begin{proof}
  (1) is immediate. For (2), one must note that the condition $[X,Y] =
  0$ implies that $X$ and $Y$ commute when $\vlie$ is regarded as a
  group by the Hausdorff series.
\end{proof}

\begin{rem}
  \label{rem:exp-remark}
  If $p \ge h$, then $\e$ defines an isomorphism $\ulie \to U$.  Since
  $\e(X)^p=\e(pX)=1$ for any $X \in \ulie$, this gives yet another
  proof that every unipotent $u \in G$ satisfies $u^p = 1$ when
  $p \ge h$.
\end{rem}

\section{Cohomology of Frobenius kernels}
\label{sec:coh-frob}

\subsection{}

Let $H$ be a linear algebraic groups over the algebraically closed
field $k$. For $d \ge 1$, denote by $H_d$ the $d$-th Frobenius kernel;
see \cite[I.9]{JRAG} for the a full discussion.  We recall some of the
details: Let $\mm \normal k[H]$ be the ideal defining 1 in the group
$H$, and put $\mm_d = \sum_{f \in \mm} k[H]f^{p^d}$.  Then $H_d$ is
the group scheme represented by the finite dimensional $k$-algebra
$k[H_d] = k[H]/\mm_d$; see \cite[I.9.6]{JRAG}. In particular, $H_d$ is
an infinitesimal group scheme \cite[I.9.6(1)]{JRAG}.

\begin{stmt}
  There is are natural bijections  
  \begin{align*}
    \Hom_{gs}(H_d,H') \iso \Hom_{gs}(H_d,{H'}_d)
    & \iso
    \Hom_{Hopf}(k[{H'}_d],k[H_d]) \\ & \iso
    \Hom_{Hopf}(\Dist(H_d),\Dist({H_d}')),
  \end{align*}
  where $\Hom_{gs}$ refers to homomorphisms of group schemes,
  $\Hom_{Hopf}$ refers to Hopf algebra homomorphisms, and $\Dist(H_d)$
  denotes the algebra of distributions of $H_d$ as in
  \cite[I.7]{JRAG}.
\end{stmt}

\begin{proof}
  Since all our group schemes are affine, the homomorphisms between
  them may be identified with comorphisms on coordinate algebras.  The
  first two isomorphisms follow from this and the fact that for any
  homomorphism $\phi:H_d \to H'$, the comorphism $\phi^*:k[H'] \to
  k[H_d]$ vanishes on $\mm'_d$.  Since $H_d$ is infinitesimal,
  $\Dist(H_d)$ identifies with the dual Hopf algebra of $k[H_d]$ by
  \cite[I.8.4]{JRAG}; the last isomorphism follows at once.
\end{proof}

\begin{stmt}
  If $A_1$ and $A_2$ are finite dimensional Hopf algebras over $k$,
  then $\Hom_{Hopf}(A_1,A_2)$ has a natural structure of algebraic
  variety over $k$.  
\end{stmt}

\begin{proof}
  Since the $A_i$ are finite dimensional, we regard $X=\Hom(A_1,A_2)$ as a
  subset of the affine space $\Aff = \Hom_k(A_1,A_2)$ of all
  $k$-linear maps.  For each $a,b \in A_1$, the map
  $\lambda_{a,b}:\Aff \to A_2$ given by $\phi \mapsto \phi(ab) -
  \phi(a)\phi(b)$ is clearly a morphism of varieties, and the set $X_a
  \subset \Aff$ of all algebra homomorphisms is the intersection of
  all $\lambda_{a,b}^{-1}(0)$, hence is a closed subvariety.
  One similarly sees that the subset $X_c$ of all coalgebra
  homomorphisms is closed, and the subset $X_{ant}$ of all antipode
  preserving linear maps is closed. Then $X = X_a \cap X_c \cap
  X_{ant}$ is also closed.
\end{proof}

\begin{stmt}
  $\Hom(H_d,H') \iso \Hom(H_d,{H'}_d)$ has the structure of an 
  $H'$-variety.
\end{stmt}

\begin{proof}
  The variety structure is evident from the previous remarks.  The
  above identification is compatible with the action of $H'$ on itself by
  inner automorphisms, and on ${H'}_d$ by the adjoint representation;
  this action yields the structure of $H'$-variety.
\end{proof}

\subsection{}

For any linear algebraic group $H$ over $k$, consider the $H$-variety
\begin{equation*}
  \A(d,H) = \Hom(\G_{a,d},H)
\end{equation*}
of the previous section, where $\G_{a,d}$ is the $d$-th Frobenius
kernel of the additive group. This is the reduced variety
corresponding to a certain (possibly not reduced) affine $k$-scheme
$\underline{\A}(d,H)$ appearing in \cite[Theorem 1.5]{SFB:infin-1-psg}
(where it is called $V_d(H)$).

\subsection{}
\label{sub:distributions}

Let $T$ be a $k$-torus with character group $X = X^*(T)$ and
co-character group $Y = X_*(T)$. Then $T$ is obtained by base change
from the $\Z$-torus $T_\Z$ defined by $\Z[X]$.  

We now use the results of \cite[I.7.8]{JRAG} to describe the algebra
of distributions. The $\Z$-algebra $\Dist(T_\Z)$ is a free
$\Z$-module; any $\Z$-basis $H_1,\dots,H_n$ of $Y$ yields a
corresponding $\Z$-basis of $\Dist(T_\Z)$: namely, all products
$\prod_{i=1}^n\dbinom{H_i}{n_i}$ with $n_i \in \N$. We will
say that the degree of such a product is $\sum_i n_i$.  The
distributions of $T$ arise by base change:
$\Dist(T) = \Dist(T_\Z) \tensor_\Z k$.

Distributions in $\Dist(T_\Z)$ are certain linear forms in
$\Hom_\Z(\Z[X],\Z)$: for $H \in Y$, $n \in \N$ and $\lambda \in X$, we
have by definition $\dbinom{H}{n}(\lambda) = \dbinom{\langle \lambda,H
  \rangle}{n}$.

Let now $T'$ (with groups $X'$, $Y'$, etc) be a second $k$-torus, and
suppose that $\phi:T \to T'$ is a morphism.  The morphism $\phi$
induces maps on the character and co-character groups: $\phi^*:X' \to
X$ and $\phi_*:Y \to Y'$.  In turn, $\phi^*$ determines a map $\Z[X']
\to \Z[X]$ and hence a morphism $\phi_\Z:T_\Z \to {T'}_\Z$ from which
$\phi$ arises by base change.

Assume that $(*)$ $\phi^*:X' \to X$ is injective and has cokernel a
finite group of order prime to $p$. This guarantees that $\dim T =
\dim T'$, $\phi$ is separable, and $\ker \phi$ is a reduced group
scheme.

The map $\Dist(\phi_\Z):\Dist(T_\Z) \to \Dist(T'_\Z)$
may be understood as follows: for $\lambda' \in X'$ we have by
definition $\Dist(\phi_\Z)\left (\dbinom{H}{n} \right )(\lambda') =
\dbinom{H}{n}(\phi^* \lambda') = \dbinom{\langle \phi^* \lambda',H
  \rangle}{n} = \dbinom{\langle \lambda',\phi_*H \rangle}{n}$.  Thus,
$\Dist(\phi_\Z)\left (\dbinom{H}{n} \right ) = \dbinom{\phi_*H}{n}$.

\begin{lem}
  Under the assumption $(*)$, $\Dist(\phi):\Dist(T) \to \Dist(T')$ is
  an isomorphism.
\end{lem}

\begin{proof} 
  Condition $(*)$ yields $\Z$-bases $H_1,\dots,H_n$ of $Y$ and
  $H_1',\dots,H_n'$ of $Y'$ and integers $a_1,\dots,a_n$ for which
  $\phi_*(H_i) = a_iH_i'$ and $\prod_i a_i \not \congruent 0
  \pmod{p}$.
  
  These bases of $Y$ and $Y'$ determine bases for the respective
  distribution algebras, and we have $\Dist(\phi_\Z)\left(\prod_i
    \dbinom{H_i}{m_i}\right) = \prod_i \dbinom{a_iH_i'}{m_i}.$ One may
  check that
  \begin{equation*}
    \prod_i
  \dbinom{a_iH_i'}{m_i} = \prod_i a_i^{m_i} \prod_i \dbinom{H_i'}{m_i} + 
    \mathcal{E},
  \end{equation*}
  where $\mathcal{E}$ is a $\Z$-linear combinations of basis elements
  of lower degree.  It follows that $\Dist(\phi) = \Dist(\phi_\Z)
  \tensor 1_k$ is an isomorphism, as claimed.
\end{proof}

\begin{theorem}
  Let $\phi:G \to G'$ be a central isogeny of connected, semisimple
  groups over $k$, as in \cite[Prop. II.1.14]{JRAG}.  Suppose that
  $\ker \phi$ is reduced. Then $\phi$ induces an isomorphism $\Dist(G)
  \iso \Dist(G')$.
\end{theorem}

\begin{proof}
  According to \cite[II.1.12(2)]{JRAG}, multiplication is an
  isomorphism 
  \begin{equation*}
    \Dist(U^-) \tensor \Dist(T) \tensor \Dist(U) \iso \Dist(G),
  \end{equation*}
  where $U$ and $U^-$ are the unipotent radicals of opposite Borel
  subgroups and $T$ is a maximal torus.  Moreover, (see
  \cite[II.1.14]{JRAG}) $\phi$ induces maps on these tensor factors;
  it is clear from the description of $\phi$ that it induces an
  isomorphism $\Dist(U) \to \Dist(U')$ (with a similar statement for
  $U^-$). Thus, it suffices to show that $\phi$ induces an isomorphism
  $\Dist(T) \to \Dist(T')$.
  
  Since $G$ and $G'$ are semisimple, $\dim T = \dim T'$; since
  $X^*(T)_\Q$ and $X^*(T')_\Q$ are spanned over $\Q$ by the roots, the
  map $\phi^*$ on character groups induced by the homomorphism
  $\phi_{\mid T}:T \to T'$ is injective; since $\coker \phi$ is reduced,
  $\ker \phi^*$ has order prime to $p$. Thus, the lemma shows that
  $\Dist(\phi_{\mid T})$ induces an isomorphism $\Dist(T) \to
  \Dist(T')$, and the result follows.
\end{proof}

The fundamental group of a root system $R$ is the finite group
$X_{sc}/\Z R$, where $X_{sc}$ is the $\Z$-lattice with basis the
fundamental dominant weights. The theorem has the following consequence:

\begin{cor}
  Let $G$ be a connected, semisimple group, with root system $R$.
  Denote the simply connected cover by $G_{sc}\to G$. If $p$ does not
  divide the order of the fundamental group of $R$, then $\A(d,G_{sc})
  \iso \A(d,G)$ as $G_{sc}$-varieties for each $d \ge 1$.
\end{cor}

\subsection{}
\label{sub:SFB-lemma}
Let $H$ a linear algebraic group over $k$ defined over $\F_p$, with
Lie algebra $\hlie$.  Let $\NN_p(\hlie)$ denote the variety of
$p$-nilpotent elements in $\hlie$.  For $d \ge 1$, put
\begin{equation}
  \label{eq:nn}
  \NN_p(d,\hlie) = \{(X_0,\dots,X_{d-1}) \mid X_i \in \NN_p(\hlie), \
  [X_i,X_j] = 0 \text{ for } 0 \le i,j < d\}.
\end{equation}
We regard $\NN_p(d,\hlie)$ as an $H$-variety with the following
action:
\begin{equation*}
  h.(X_0,X_1,\dots,X_{d-1}) = (\Ad(h)X_0,\Ad(Fh)X_1,\dots,
  \Ad(F^{d-1}h)X_{d-1}),
\end{equation*}
where $F$ denotes the Frobenius morphism on $H$.

We have the following analogue of Theorem \ref{sub:distributions}.
\begin{lem}
  Let $G$ be connected, semisimple with root system $R$ and simply
  connected cover $G_{sc}$. If $p$ does not divide the order of the
  fundamental group of $R$, there is an isomorphism of
  $G_{sc}$-varieties $\NN_p(d,\Lie(G)) \iso \NN_p(d,\Lie(G_{sc})$ for
  each $d \ge 1$.
\end{lem}

\begin{proof}
  This follows from the observations made in
  \cite[0.13]{hum-conjugacy}.
\end{proof}

The following result was obtained in \cite[Lemma 1.7]{SFB:infin-1-psg}.
\begin{prop}
  Suppose that $H$ has a faithful rational representation $(\rho,V)$
  with the property that $\exp(d\rho(X)) \in H$ for each $X \in
  \NN_p(\hlie)$, where $\exp(d\rho(X))$ is the (truncated) exponential
  in $\GL(V)$.  Then there is an isomorphism of $H$-varieties
  $\NN_p(d,\hlie) \iso \A(d,H)$.
\end{prop}
 
Actually, in \emph{loc. cit.}, one gets an isomorphism of schemes
$\underline{\NN}_p(d,\hlie) \iso \underline{\A}(d,H)$; for each
commutative $k$-algebra $\Lambda$,
$\underline{\NN}_p(d,\hlie)(\Lambda)$ is the set described by
\eqref{eq:nn} except that the $X_i$ are taken from $\hlie\tensor_k
\Lambda$.  To get this isomorphism of schemes, one must make the
assumption that the exponential of any $p$-nilpotent $X \in \hlie
\tensor_k \Lambda$ lies in $H(\Lambda)$ for each $\Lambda$.  A look at
the proof in \cite{SFB:infin-1-psg} shows that we still get an
isomorphism of varieties with our weaker assumption.

If $(\rho,V)$ satisfies the hypothesis of the lemma, we 
say that it is an exponential-type representation of $H$.

\subsection{}
\label{sub:scheme-iso}
Let now $G$ be a connected, semisimple, algebraic group over $k$.
Let $\NN_p = \NN_p(\glie)$, $\NN_p(d) = \NN_p(d,\glie)$, and 
$\A(d) = \A(d,G)$.

The group $G$ is determined up to isogeny by its root system $R$.  Let
$r$ denote the rank of $G$. When $G$ is quasisimple, $R$ is either of
classical type (hence is one of $A_r$, $B_r$, $C_r$, or $D_r$), or $R$
is of exceptional type (hence is one of $E_6$, $E_7$, $E_8$, $F_4$ or
$G_2$).

\begin{theorem} Let $G$ be semisimple. Then $\NN_p(d) \iso \A(d)$ as
  $G$-varieties in each of the following cases:
  \begin{enumerate}
  \item $p > 2h-2$ where $h$ is the Coxeter number.
  \item $G$ is quasisimple and $R$ is of classical type, and
    moreover $p \not = 2$ if $R = B_r,C_r,D_r$, and $r \not
    \congruent -1 \pmod{p}$ if $R = A_r$.  
  \item $G$ is quasisimple and  $R$ is of exceptional
    type, and moreover $p \ge p_0$ where $p_0$ is given in the
    following table:
    \begin{equation*}
      \begin{array}[t]{ll}
        R & p_0 \\
        \hline \hline
        G_2 & 7 \\
        F_4 & 17 \\
        \end{array} \quad
      \begin{array}[t]{ll}
        R & p_0 \\
        \hline \hline
        E_6 & 17 \\
        E_7 & 29 \\
        E_8 & 59 \\
      \end{array}
    \end{equation*}
  \end{enumerate}
\end{theorem}

\begin{proof}
  In all three cases, the condition on $p$ guarantees that it does not
  divide the order of the fundamental group (this is well-known, and
  may be checked by looking at the tables in \cite{BouLie456}).  So it
  suffices by Corollary \ref{sub:distributions}, Lemma
  \ref{sub:SFB-lemma} and Proposition \ref{sub:SFB-lemma} to show that
  there is a semisimple group $G'$ isogenous to $G$, and an
  exponential-type representation $(\rho,V)$ of $G'$.
  
  In case (1), this follows from Corollary \ref{sub:admiss-group}
  together with Remark \ref{sub:admiss-group}(1).
  
  In case (2), this follows from \cite[Lemma 1.8]{SFB:infin-1-psg}.

  In case (3), we get the result again by Corollary \ref{sub:admiss-group}
  together with Remark \ref{sub:admiss-group}(2).
\end{proof}

\subsection{}
\label{sub:intrinsic}
Let $G$ be connected and reductive, and let $P \le G$ be a parabolic
subgroup with unipotent radical $\V \le P$ and $\vlie = \Lie(\V)$.
Let $n(P)$ be the integer defined as in \ref{sub:nilp-class}; this
coincides with the nilpotence class of $\V$ and $\vlie$ provided $p$
is good; in particular, it coincides with the nilpotence class of
$\V_\Q$ and $\vlie_\Q$ as in section \ref{sec:exponential-big-p}. The
following is related to the question posed in \cite[Remark
1.9]{SFB:infin-1-psg} (see the discussion in the introduction).

\begin{theorem}
  Assume that $n(P) < p$.
  Then there is an injective morphism of $P$-varieties
  \begin{equation*}
    \e:\NN_p(d,\vlie) \to \A(d,\V).
  \end{equation*}
\end{theorem}

\begin{rem}
  The point of the theorem is that $\e$ does not depend on the choice
  of a faithful representation of $G$.
\end{rem}

\begin{proof}
  In this case, we must first work with schemes in order to know that
  the map we define is a morphism.

  Recall from Corollary \ref{cor:exp-iso} that there is an isomorphism
  of group schemes $\e:\vlie \to \V$.  Thus for each $k$-algebra
  $\Lambda$, each $X \in \vlie \tensor_k \Lambda$ determines a
  homomorphism $\e_X:\G_{a,\Lambda} \to \V_\Lambda$. If
  $\vec{X}=(X_0,X_1,\dots,X_{d-1}) \in
  \underline{\NN}(d,\vlie)(\Lambda)$, one emulates the construction in
  \cite[Remark 1.3]{SFB:infin-1-psg} to obtain a homomorphism of group
  schemes $\e_{\vec{X}}:\G_{a,\Lambda} \to \V_\Lambda$ (note that we
  must use here (2) of Corollary \ref{cor:exp-iso}), and hence (by
  ``restriction'') a homomorphism of group schemes
  $\e_{\vec{X}}:\G_{a,d,\Lambda} \to \V_\Lambda$.
  
  It is now easy to see that the assignment $\vec{X} \mapsto
  \e_{\vec{X}}$ is functorial in $\Lambda$, hence defines a morphism
  of schemes $\underline{\NN}(d,\vlie) \to \underline{\A}(d,\V)$. We
  get then also a morphism of varieties $\NN(d,\vlie) \to \A(d,\V)$;
  $P$-equivariance follows from (1) of Corollary \ref{cor:exp-iso}.
  
  To prove injectivity, we essentially copy the proof of \cite[Lemma
  1.7]{SFB:infin-1-psg}. Suppose that $\e_{\vec{X}} = \e_{\vec{Y}}$.
  Differentiating gives then $X_0 = Y_0$. Multiplying each
  homomorphism with $\e_{-X_0}$, one sees that
  $\e_{(0,X_1,\dots,X_{d-1})} = \e_{(0,Y_1,\dots,Y_{d-1})}$ are equal.
  But then $\e_{(X_1,\dots,X_{d-1})}$ and $\e_{(Y_1,\dots,Y_{d-1})}$
  coincide in $\A(d-1,\V)$, and the injectivity of $\e$ follows by
  induction (note that $\e$ is an isomorphism of varieties when $d=1$;
  see \cite[Lemma 1.6]{SFB:infin-1-psg}).
\end{proof}

\subsection{}
Let $H$ be a linear algebraic group over $k$.  We recall briefly the
significance of the variety $\A(d,H)$ for the cohomology of $H_d$.  In
the papers \cite{SFB:infin-1-psg} and \cite{SFB:support-var}, Suslin,
Friedlander and Bendel define a ring homomorphism
\begin{equation*}
  \Phi:H^{\text{even}}(H_d,k) \to k[\underline{\A}(d,H)],
\end{equation*}
and they show that the map induced by $\Phi$ on the corresponding
schemes is a topological homeomorphism; clearly
the same is still true after replacing $\underline{\A}(d,H)$
with $\A(d,H)$.

Now consider a connected reductive group $G$ over $k$.  It was shown
by Friedlander and Parshall that for sufficiently large $p$, the
cohomology $H^i(G_1,k)$ vanishes for $i$ odd, and that the even
cohomology ring $H^{\text{even}}(G_1,k)$ may be identified with the
graded coordinate ring of $\NN(\glie)$.  See also \cite{AJ1}, where it
is proved that $p > h$ is sufficient. The proof of this fact in
\emph{loc. cit.} relies on knowledge of the dimensions of the
homogeneous parts of $k[\NN(\glie)]$; thus it seems likely that
further understanding of the cohomology of $G_d$ might benefit from
some understanding of $\NN_p(d,\glie)$. We conclude the paper with the
following, in the hope that it might be useful (compare
\cite[3.9]{AJ1}).

\begin{prop}
  If $p$ is a good prime for $G$, there is an injective
  homomorphism of $G$-modules
    \begin{equation*}
      k[\NN_p(d,\glie)] \to H^0(G/B,k[\NN_p(d,\ulie)]).
    \end{equation*}
\end{prop}

\begin{proof}
  We mimic the argument in \cite{AJ1}. Let $X = \bigoplus_{i=0}^{d-1}
  \glie^{[i]}$ (where the exponent $[i]$ denotes the $i$-th Frobenius
  twist); $G$ acts on $X$ by $\alpha = \bigoplus \Ad^{[i]}$.  We
  denote also by $\alpha$ the action of $G$ on the algebra $k[X]$ of
  regular functions on $X$.  There is a homomorphism
  \begin{equation*}
    k[X] 
    \to H^0(G/B,k[\NN_p(d,\ulie)])
  \end{equation*}
  obtained by mapping $f \in k[X]$ to the section $g \mapsto
  \left(\alpha(g^{-1}) f\right)_{\mid \NN_p(d,\ulie)}$.  The kernel is
  \begin{equation*}
    \{f \in k[X] \mid f \text{\ \ vanishes on $\alpha(g)\NN_p(d,\ulie)$
      for all $g \in G$}\},
  \end{equation*} 
  and we claim this is the vanishing ideal of $\NN_p(d,\glie)$.  It
  suffices to see that if $\{X_i\}$ is a set of pairwise commuting
  nilpotent elements in $\glie$, then the Abelian Lie algebra $\alie$
  which they span is contained in a Borel subalgebra; that is a
  consequence of the lemma which follows.
\end{proof}

\begin{lem}
  Suppose that $p$ is good for $G$, and that $G$ is a Borel subgroup
  of $G$ with unipotent radical $U$. Let $\alie \subset \glie$ be an
  Abelian subalgebra generated by nilpotent elements.  Then there is a
  $g \in G$ such that $\Ad(g)\alie \subset \ulie = \Lie(U)$.
\end{lem}

\begin{proof}
  There is a central isogeny $G' \to G$ where $G'$ is a direct product
  of a torus and quasisimple groups satisfying the hypothesis of
  Proposition \ref{sub:spaltenstein}. This isogeny induces a bijection
  (not in general an isomorphism of varieties) on the nilpotent sets
  in the respective Lie algebras. Thus, we may replace $G$ with $G'$,
  so that we may apply the results of \cite{spalt-transverse}.

  The result is well known if $\dim \alie = 1$. So now suppose that
  $\dim \alie > 1$, and let $0 \not = X \in \alie$.  By the Theorem
  proved in \cite{spalt-transverse}, there is a proper parabolic
  subgroup $P$ of $G$ with Levi decomposition $P = \Levi\V$ such that
  $X \in \vlie = \Lie(\V)$ and $\clie_\glie(X) \subset \plie =
  \Lie(P)$ [in general, $X$ need not be a Richardson element in
  $\vlie$].  Thus, we have $\alie \subset \clie_\glie(X) \subset
  \plie$.  Since $X \in \vlie$, the image of $\alie$ in $\plie /
  \vlie$ has dimension strictly less than that of $\alie$. We obtain
  by induction on $\dim \alie$ some $g \in \Levi$ such that the image
  of $\alie$ in $\llie \iso \plie / \nlie_\plie$ is conjugate via
  $\Ad(g)$ to a subalgebra of a Borel subalgebra of $\llie$.  Since
  $\Levi$ leaves $\vlie$ invariant, $\Ad(g)(\alie)$ is contained in a
  Borel subalgebra of $\plie$ (which is in turn a Borel subalgebra of
  $\glie$).  Since all Borel subalgebras of $\glie$ are conjugate, we
  obtain the lemma.
\end{proof}

\providecommand{\bysame}{\leavevmode\hbox to3em{\hrulefill}\thinspace}

\end{document}